\documentclass[12pt]{article}
\usepackage{fullpage}
\usepackage{amssymb}
\usepackage{latexsym}
\newcommand{\R}{\mathbb{R}}
\newcommand{\tz}{t_0}
\newcommand{\splitstat}{(\alpha^{-1} w_0^+, - w_0^-)}
\newcommand{\lplp}{L^p(\Omega) \times L^p(\Omega)}
\newcommand{\conel}{C^{1, \lambda}(\overline{\Omega})}
\newcommand{\RN}{\R^N}
\newcommand{\N}{\mathbb{N}}
\newcommand{\coneob}{C^1(\overline{\Omega})}
\newcommand{\wsi}{\overline{w}^i}
\newcommand{\xt}{(x,t)}
\newcommand{\pwk}{(P_{w^k})}
\newcommand{\ukvk}{(\uk, \vk)}
\newcommand{\uk}{u^k}
\newcommand{\pks}{(P_k(S))}
\newcommand{\ukj}{u^{k_j}}
\newcommand{\fit}{\tilde{f}_1(u,v,t)}
\newcommand{\ftt}{\tilde{f}_2(u,v,t)}
\newcommand{\bc}{\alpha m_1 - m_2}
\newcommand{\jz}{j_0}
\newcommand{\conez}{C^1_0(\overline{\Omega}) \times C^1_0(\overline{\Omega})}
\newcommand{\ucomb}{u- \alpha^{-1}v}
\newcommand{\done}{\delta_1}
\newcommand{\ucombp}{(\ucomb)^+}
\newcommand{\vcomb}{v- \alpha u}
\newcommand{\vcombp}{(\vcomb)^+}
\newcommand{\hkj}{h^{k_j}}
\newcommand{\lp}{L^p(\Omega)}
\newcommand{\lkj}{l^{k_j}}
\newcommand{\hk}{h^k}
\newcommand{\wh}{\hat{w}}
\newcommand{\llk}{l^k}
\newcommand{\wz}{w_0}
\newcommand{\vk}{v^k}
\newcommand{\fixt}{(\cdot, t)}
\newcommand{\vkj}{v^{k_j}}
\newcommand{\pf}{\noindent {\em Proof.}~~}
\newcommand{\pk}{(P_k)}
\newcommand{\pkt}{(P_{k,T})}
\newcommand{\kj}{k_j}

\newcommand{\tj}{t_j}
\newcommand{\kz}{k_0}
\newcommand{\wtp}{W^{2,p}(\Omega)}
\newcommand{\wtpz}{W^{2,p}_0(\Omega)}

\newcommand{\me}{\mathcal{E}}

\newcommand{\nw}{N_{\delta}(\wz)}
\newcommand{\npair}{N_{\delta_1}(\alpha^{-1} \wz^+, - \wz^-)}

\newcommand{\tw}{\tilde{w}}
\newcommand{\xj}{x_j}
\newcommand{\eh}{\hat{\eta}}
\newcommand{\ra}{\to}
\newcommand{\brw}{{\mathcal{B}}_R(w)}
\newcommand{\ms}{{\mathcal{S}}}

\newcommand{\epz}{\epsilon_0}

\newcommand{\epj}{\epsilon_j}
\newcommand{\tb}{\bar{t}}
\newcommand{\xb}{\bar{x}}
\newcommand{\xn}{x_n}
\newcommand{\xtt}{\tilde{x}}
\newcommand{\ut}{\tilde{U}}
\newcommand{\zt}{\tilde{z}}
\newcommand{\ep}{\epsilon}
\newcommand{\oj}{\Omega_j}
\newcommand{\stat}{(S)}
\newcommand{\qed}{\hfill $\Box$ \\[5pt]}
\newcommand{\wkp}{(\wk)^+}
\newcommand{\wkm}{(\wk)^-}
\newcommand{\wk}{w^k}
\newcommand{\wkj}{w^{\kj}}
\newcommand{\bfive}{{\bf (b5)}}
\newcommand{\btwo}{{\bf (b2)}}
\newcommand{\bone}{{\bf (b1)}}
\newcommand{\bthree}{{\bf (b3)}}

\newcommand{\mkone}{m_1^k}
\newcommand{\mktwo}{m_2^k}

\newcommand{\mki}{m_i^k}
\newcommand{\mkjone}{m_1^{\kj}}
\newcommand{\mkjtwo}{m_2^{\kj}}
\newcommand{\half}{\frac{1}{2}}
\newcommand{\lk}{\Lambda^k}
\newcommand{\lj}{\Lambda^{k_j}}

\newcommand{\clo}{\overline{\Omega}}
\newcommand{\lto}{L^2(\Omega)}
\newtheorem{definition}{Definition}[section]
\newtheorem{theorem}[definition]{Theorem}

\newtheorem{lemma}[definition]{Lemma}

\begin{document}
\title{On long-time dynamics for competition-diffusion systems with inhomogeneous
Dirichlet boundary conditions}
\author{E.C.M. Crooks, E.N. Dancer and D. Hilhorst}
\date{}
\maketitle
\begin{abstract}
We consider a two-component competition-diffusion system with equal diffusion
coefficients and inhomogeneous Dirichlet boundary conditions. When the interspecific
competition parameter tends to infinity, the system solution converges to that of a 
free-boundary problem. If all stationary solutions of this limit problem are non-degenerate
and if a certain linear combination of the boundary
data does not identically vanish, then for sufficiently large interspecific competition, 
all non-negative solutions of the competition-diffusion system
converge to stationary states as time tends to infinity. 
Such dynamics are much simpler than
those found for the corresponding system with either homogeneous Neumann or homogeneous Dirichlet boundary conditions. \\[10pt]
\end{abstract}

\section{Introduction}
In this paper, we show that, under certain conditions, the competition-diffusion system
\begin{eqnarray}
u_t & = & \Delta u + f(u) - kuv ~~~\mbox{in}~\Omega, \nonumber \\
v_t & = & \Delta v + g(v) - \alpha kuv~~~ \mbox{in}~ \Omega, \label{introsyst}
\end{eqnarray}
with inhomogeneous Dirichlet boundary conditions
\begin{eqnarray}
u & = & m_1 \geq 0  \nonumber ~~~\mbox{on}~\partial \Omega,\\
v & = & m_2 \geq 0 \label{introbc}~~~~\mbox{on}~ \partial \Omega,
\end{eqnarray}
has simple long-time dynamics for large positive values of the competition
parameter $k$. Here $\Omega \subset \RN$ is smooth and bounded, $f$ and $g$
are positive on $(0,1)$ and negative elsewhere, and $\alpha >0$. 
Such reaction-diffusion systems are well-known in the modelling of competition between
two species of population densities $u(x,t)$ and $v(x,t)$, and we refer to the introduction
of \cite{CDHMN} for a brief review. These models can be used to study the dynamics of
the spatial segregation between the competing species. 
The parameters $k$ and $\alpha$
may be thought of as representing the interspecific competition rate and the competitive advantage of $v$ over $u$ respectively. 
Zero flux (that is, zero Neumann) are the most commonly imposed boundary conditions.
But when the two species have quite different preferences for environmental conditions, then competition occurs mainly in a region $\Omega$
where their habitats overlap  and this gives rise to  boundary conditions
(\ref{introbc})  on $\partial \Omega$ \cite{NambaMimura}. 

More precisely, 
we prove that if $\alpha m_1 - m_2$ is not identically zero on $\partial \Omega$
and all stationary solutions of the limit problem
\begin{eqnarray}
- \Delta w & = & \alpha f (\alpha^{-1} w^+) - g(-w^-) =: h(w) ~~~\mbox{in}~\Omega, \nonumber \\
w & = & \alpha m_1 - m_2 ~~~\mbox{on}~\partial \Omega, \label{introlimit}
\end{eqnarray}
are non-degenerate (see Definition \ref{defnondegen}), then for $k$ sufficiently large, all non-negative solutions of (\ref{introsyst})
approach stationary states as $t \rightarrow \infty$. 

Two remarks on our hypotheses and results should be made at the outset. 
First, provided we suppose that $\alpha m_1 - m_2$ is not identically zero on $\partial \Omega$,
 our system  (\ref{introsyst}) with inhomogeneous Dirichlet boundary conditions has much simpler dynamics than the corresponding system with
zero Dirichlet or Neumann boundary conditions.
 In \cite{DZ}, it is observed that such systems may have solutions that are small
($O (1/k)$) for all time. To ensure that for large $k$ these solutions converge to a stationary
solution of the $k-$dependent system as $t \rightarrow \infty$, it is necessary to impose a condition
of there being no ``circuits'' of positive heteroclinic orbits of an associated limit system
(see \cite[Assumption C3]{DZ} and \cite{Dpre}) plus a condition on a linear limit problem
(\cite[Assumption C1]{DZ}). No such additional assumptions are needed here to show simple dynamics.
Compare Theorem \ref{thmsimpledynamics} with
\cite[Thm 5]{DZ}.

Second,  the condition that all solutions of the stationary
limit problem (\ref{introlimit}) are non-degenerate does not always hold for our boundary conditions - not
even in one space dimension. 
This contrasts with the case of zero Neumann or zero
 Dirichlet boundary conditions, in which non-degeneracy does hold in one space-dimension - see a remark in \cite[p 472]{DZ}.
But some genericity results can be shown for our inhomogeneous Dirichlet case, and we discuss
these, together with the possible failure of non-degeneracy in one dimension, in Section 6.

Our methods owe much to \cite{DZ}, which treats (\ref{introsyst}) with zero Neumann
boundary conditions.
The idea is 
first to use a blow-up method to show that for each $\delta >0$, one of $u$ or/and $v$ must be
small at each $(x,t) \in \Omega \times [\delta, \infty)$ for sufficiently large $k$ (Section 2). 
This results in the linear combination $w=\alpha u - v$ satisfying the scalar
equation
\begin{eqnarray}
 w_t & = & \Delta w + h(w) + \mathcal{O}_k(1), ~~~\mbox{in}~\Omega, \label{introlya}\\
 w & = & \alpha m_1 - m_2 ~~~\mbox{on}~ \partial \Omega \nonumber,
 \end{eqnarray}
 where $\| \mathcal{O}_k(1)\|_{L^2(\Omega)} \rightarrow 0$ as $k \rightarrow \infty$ uniformly
 in $t \in [\delta, \infty)$. Note that here we can only estimate the $L^2$-norm of $\mathcal{O}_k(1)$, rather than the $L^{\infty}$-norm, as in \cite{DZ}, if the given boundary data $m_1, m_2$ is not assumed to be segregated on $\partial \Omega$. But this $L^2$-estimate is sufficient to study the long-time behaviour of (\ref{introsyst}). The Lyapunov function for (\ref{introlya}) with $\mathcal{O}_k(1) =0$ can then be used (Section 3) to show that $w$ must lie close to solutions of
 (\ref{introlimit}) for $k, t$ large, under the condition that solutions of (\ref{introlimit}) are
 isolated in $L^2(\Omega)$. Section 4 then shows that if these stationary solutions are in fact
 all non-degenerate, then solutions of (\ref{introsyst}) must approach stationary states
 of (\ref{introsyst}) as $t \rightarrow \infty$. 
  Note
  that the non-degeneracy required in Section 4 does imply the isolatedness used in Section 3, even though the function
 $h$ in (\ref{introlimit}) being only locally Lipschitz at its zero set means that the inverse function theorem cannot be applied directly to the operator $w \mapsto \Delta w + h(w)$ (see, for example, remark {\em (ii)} at the end of Section 6).
 That there is a (locally) unique stationary solution
 of (\ref{introsyst}) close to $(\alpha^{-1} w^+, -w^-)$ for $w$ a non-degenerate solution of
 (\ref{introlimit}) is shown in Section 5 using index-theory arguments similar to those
 in \cite{DG, DD}. Our inhomogeneous boundary values here necessitate careful modification of
 various arguments in \cite{DZ, DG, DD}, particularly the blow-up argument in Section 2, and
 also in the bounds and index arguments used to prove local uniqueness in Section 5. 
 Section 6 is devoted to non-degeneracy of stationary solutions of (\ref{introlimit}), as mentioned
 above. We use an approach 
 from \cite{ST, Dgen} to show that all stationary solutions of the limit problem are non-degenerate
 for generic boundary data by applying the version of Sard's Theorem from \cite{Q} to a
 suitable map. Our function $h$ defined in (\ref{introlimit}) is locally Lipschitz but
 not in general continuously differentiable; \cite{Dgen} extends the work of \cite{ST} to deal with
 such non-smooth functions, and we use the ideas from \cite{Dgen} here.  

This paper follows on from the related work \cite{CDHMN}, in which a spatial segregation
limit is derived for the generalisation of (\ref{introsyst}) in which the diffusion coefficients
of $u$ and $v$ are allowed to differ. It is shown there  that for each $T>0$, $u$ and $v$ converge in $L^2(\Omega \times (0, T))$ as $k \rightarrow \infty$, where in the limit, $uv=0$ almost everywhere and $w = \alpha u - v$ is the solution of a limiting free boundary problem.
Here, our assumption that the diffusion coefficients of $u, v$ are in fact the same enables
us to form the equation (\ref{introlya}) which plays a key r\^ole in the rest of our analysis. 
It also allows
us to establish the key estimates in Section 1 {\em uniformly} in $t$, which enables us to use
the Lyapunov-function argument in Section 3.  (Note that the argument in Section 1 yields
estimates uniform in the unbounded time-interval $t \in [\delta, \infty)$ for each $\delta >0$ 
and we exploit this in the energy argument that is given in Section 3. But estimates uniform
on bounded time intervals would in fact be sufficient to obtain the result in Section 3 using
a slightly different argument exploiting \cite[Theorem 3.4.1]{Henry} - see \cite{DZ}.)
\\[0.5cm]
\noindent {\bf Acknowledgements} Elaine Crooks gratefully acknowledges financial support from the Michael Zilkha Trust, of Lincoln College, Oxford. Norman Dancer was partially supported by the Australian Research Council.  Danielle Hilhorst was partially supported by the RTN contract
FRONTS-SINGULARITIES HPRN-CT-2002-00274.

\section{Formulation of the problem and a key lemma}
\label{one} Let $\Omega$ be a bounded, open, connected subset of
$\RN$ with  boundary $\partial \Omega$ of class $C^{2, \mu}$ for
some $\mu > 0$
and $Q := \Omega \times \R^+$. Let $k \in \N$ and consider the $k-$dependent problem\\[10pt]
$$(P_k) \left\{
\begin{array}{ll}
u_t = \Delta u + f(u) - kuv & \mbox{in}~Q, \\
v_t = \Delta v + g(v) - \alpha kuv & \mbox{in}~Q, \\
u=m_1^k & \mbox{on}~\partial \Omega \times \R^+,\\
v=m_2^k &  \mbox{on}~\partial \Omega \times \R^+,\\
u(x,0)= \uk_0(x), ~v(x,0) = \vk_0(x) & \mbox{for}~x \in \Omega,
\end{array} \right. $$\\[5pt]
where it is supposed throughout that
\begin{description}
\item[(a)] $f$ and $g$ are continuously differentiable functions on $[0, \infty)$
such that $f(0)=g(0)=0$ and $f(s)<0$, $g(s)<0$ for all $s>1$;
\item[(b1)] $\mkone, \mktwo \geq 0$ and $\mkone, \mktwo \in  \wtp$ where $p > N$;
\item[(b2)] $\mkone, \mktwo$ are bounded in $\wtp$ independently of $k$; 
\item[(b3)] there exist $m_1, m_2 \in \wtp$ such that $\alpha m_1 - m_2$ is not identically zero
on $\partial \Omega$ and
   $$\mkone \rightarrow m_1~~~\mbox{and}~~~ \mktwo \rightarrow m_2 ~~~~\mbox{in}~~~
   C^{1, \lambda^{\prime}} (\overline{\Omega})$$
 for each  $\lambda^{\prime} \in (0, \lambda)$, with $\lambda := 1 - N/p$ ({\em cf.} \cite[p 47]{Temam});
\item[(b4)] the initial conditions $u^k_0$ and $v_0^k$
are defined by $$u_0^k(x) = \mkone (x),~~~~~v_0^k(x) = \mktwo (x)~~~\mbox{for}~~x \in \Omega.$$

\end{description}
Some of our results will need the following stronger hypothesis on the limiting
boundary behaviour of $(u,v)$;
\begin{description}
\item[(b5)]  let $\Gamma_1, \Gamma_2$ be closed smooth sub-manifolds-with-boundary
of $\partial \Omega$, with non-empty interior in
$\partial \Omega$, and such that $\partial \Omega = \Gamma_1 \cup \Gamma_2$. 
Then  $m_1$ and $m_2$ in {\bf(b3)}
 are such that $m_i \equiv 0$ on $\Gamma_j$ where $j \neq i$.
\end{description}
\smallskip
We note the following basic consequence of  \btwo~ and \bthree.

\begin{lemma}
\label{lemA}
Suppose $\Omega \ni x_k \rightarrow x$ as $k \rightarrow \infty$.  Then $m_i^k(x_k) \rightarrow m_i(x)$ as $k \rightarrow \infty$
for $i \in \{1,2\}$. 
\end{lemma}
\pf \btwo~ implies that given $\epsilon >0$, there exists $k_0 >0$ such that $|m_i^k(x_k) - m_i(x_k)| < \epsilon /2$
for all $k \geq k_0$. And since $m_i$ is continuous, by \bthree, $|m_i(x_k) - m_i(x)| < \epsilon /2$ for all $k$
sufficiently large. The result follows. \qed

\noindent By a solution of problem $\pk$ we will mean a pair $(u,v)$ such that
$u,v \in C(\overline{Q})\cap C^{2,1}(\Omega \times [\tz, \infty))$ for any $\tz >0$.
We will say that $(u,v)$ is a solution of problem $(P_{k,T})$ if
$u,v \in C(\overline{Q})\cap C^{2,1}(\Omega \times [\tz, T])$ for $\tz \in (0,T)$
satisfies $\pk$ with $\R^+$ replaced by $(0,T)$.
\\[10pt]
We begin with some standard preliminaries on {\em a priori} bounds and global
well-posedness for the problem $\pk$.
\begin{lemma}
\label{apriori}
Let $M \geq \max \{1, \mkone, \mktwo \}$ and suppose that
$(\uk, \vk)$ is a solution of $(P_{k,T})$ for some $T>0$. Then
$$ 0 \leq \uk, \vk \leq M ~ ~ ~  \mbox{in} ~ ~ \overline{\Omega} \times [0,T].$$
\end{lemma}

\pf Define $\mathcal{L}_1(u)= u_t - \Delta u -f(u)+ kuv$ and $\mathcal{L}_2(v)=
v_t - \Delta v - g(v)+\alpha kuv$. Since $\mathcal{L}_i(0)=0$, 
 $i=1,2$, it follows from the maximum principle that $\uk, \vk \geq 0$. One can then
 check that $\mathcal{L}_i(M) \geq 0$, $i = 1,2$, which completes the proof of
 Lemma \ref{apriori}.
\qed

\begin{lemma}
\label{global}
There exists a unique solution $\ukvk$ of $\pk$ for each $k \in \N$.
\end{lemma}

\pf By \cite[Thms 9.15 and 9.19]{GT}, there exist $h_1, h_2 \in
C^{\infty}(\Omega) \cap \wtp$ such that for $i \in {1,2}, \Delta
h_i =0$ in $\Omega$ and $h_i = \mki $ on $\partial \Omega$ (in the
sense of trace).
Defining $U:= u-h_1$, $V:= v-h_2$ allows application of \cite[Prop 7.3.2]{Lunardi} to
the corresponding system for $U$ and $V$ with homogeneous boundary conditions to
yield the existence of a unique solution $\ukvk$ of $\pkt$ for some $T>0$. That we can take
$T=\infty$ follows from the {\em a priori} bounds of Lemma \ref{apriori} and the last
part of \cite[Prop 7.3.2]{Lunardi}. \qed

\noindent Given the solution $\ukvk$ of $\pk$, define
\begin{equation}
\label{wk}
\wk = \alpha \uk - \vk.
\end{equation}
Then $\wk$ satisfies the equation
\begin{equation}
\label{wkeqn}
\pwk \left\{
\begin{array}{ll}
\wk_t = \Delta \wk + \alpha  f(\uk) - g(\vk) & \mbox{in}~~Q,\\
\wk = \alpha \mkone - \mktwo & \mbox{on}~~\partial \Omega \times \R^+, \\
\wk (x,0) = \alpha \uk_0 (x) - \vk_0 (x) & \mbox{for}~~x \in \Omega.
\end{array}
\right.
\end{equation}
Note that the explicitly $k-$dependent terms in $\pk$ cancel on forming the equation for
$\wk$. Together with Lemma \ref{apriori}  and \btwo,~ this gives $k$-{\em independent} bounds
for $\wk$ which are crucial in the following. \\[10pt]

\noindent Now fix $\beta >0$ and $\xi \in (0, \half)$ and define
\begin{equation}
\label{komega}
\lk:= \left\{ x \in \Omega : \mbox{dist}(x, \partial \Omega) \geq \frac{\beta}{k^{\half - \xi}} \right\}.
\end{equation}
The following lemma is crucial.

\begin{lemma}
\label{keylem} $(i)$ 
 Let $\epsilon, M, \tz >0$. Then there exists $k_0 > 0$ such that if $k
  \geq k_0$ and $\ukvk$ is a solution of $\pk$ on $\Omega \times (0, \infty)$ with $0 \leq
  \uk,\vk \leq M$, then given any $x \in \lk$ and $t \geq \tz$, either
\begin{equation}
\label{small}
\uk(x,t) \leq \epsilon~~~\mbox{ or}~~~ \vk(x,t)\leq \epsilon.
\end{equation}
$(ii)$ If, in addition,  $\mkone, \mktwo$ satisfy the supplementary condition \bfive, then (\ref{small})  holds
for any $x \in \overline{\Omega}$ and $t \geq \tz$.
\end{lemma}

\pf We adapt the blow-up argument used in the proof of \cite[Thm 1]{DZ}, and will prove parts
$(i)$ and $(ii)$ in parallel. The first step is to use a contradiction argument to obtain a
limiting system for part $(i)$ that is defined on $\RN \times \R$, and then to use a similar argument
for part $(ii)$ to obtain the same limiting equations  but  defined on either $\RN \times \R$ or on
$H \times \R$ for a half-space $H$. The second step will be to show that for both possible limit
problems, one component must vanish identically, which will lead to a contradiction.  \\[6pt]
\noindent
First consider part $(i)$ and
suppose, for contradiction, that there exist $\epsilon_0>0$, $j-$indexed
sequences  $\kj \to \infty, ~\tj \geq \tz$, $\xj \in \lj$ and solutions $u^{\kj},
v^{\kj}$ of $(P_{\kj})$ such that $u^{\kj}(\xj, \tj) \geq \epsilon_0$ and  $v^{\kj}(\xj, \tj) \geq \epsilon_0$.\\[6pt]
\noindent
Define new $j-$dependent variables $x'=\sqrt{\kj}(x-\xj), t'= \kj (t-\tj)$ and
the  sets $\Omega_j$ by $x' \in \Omega_j$ whenever $x \in \Omega$.
Then for $x' \in \Omega_j$ and $t' \in [- \kj \tj, \infty)$, the functions $U^{\kj},
V^{\kj}$ defined by
$$(U^{\kj}, V^{\kj})(x',t')= (U^{\kj}, V^{\kj} )(\sqrt{\kj}(x-\xj),\kj (t-\tj)) =
(u^{\kj}, v^{\kj})(x,t)$$
satisfy

\begin{equation} \left\{
\begin{array}{ll}
U^{\kj}_t = \Delta U^{\kj} +\kj^{-1} f(U^{\kj}) - U^{\kj}V^{\kj} & \mbox{in}~\Omega_j
\times [- \kj \tj, \infty) , \\
V^{\kj}_t = \Delta V^{\kj} +\kj^{-1} g(V^{\kj}) - \alpha U^{\kj}V^{\kj} & \mbox{in}~\Omega_j \times  [- \kj \tj, \infty), \\
U^{\kj}=\mkjone & \mbox{on}~\partial \Omega_j \times [- \kj \tj, \infty) ,\\
V^{\kj}=\mkjtwo &  \mbox{on}~\partial \Omega_j \times [- \kj \tj, \infty),\\
U^{\kj}(x',- \kj \tj)= u^{k_j}_0(x), ~V^{\kj}(x',- \kj \tj) = v^{k_j}_0(x) & \mbox{for}~x' \in \Omega_j.
\end{array} \right. \label{nveq}
\end{equation}
Note that $0 \in\Omega_j$, $U^{\kj}(0,0)\geq \epsilon_0$ and $V^{\kj}(0,0)\geq \epsilon_0$.\\[6pt]
\noindent
Consider what happens to the system (\ref{nveq}) as $j \to \infty$. Note first that since
$t_j \geq \tz$ for each $j$ and $k_j \to \infty$,  $[- \kj \tj, \infty)$ tends to
$\R$ as $j \to \infty$, in the sense that given a compact interval $I \subset \R$, there
exists $j_0$ such that $I \subset [ -\kj \tj, \infty)$ for all $j \geq j_0$.
This will enable us to obtain limiting problems defined
for {\em all} $t \in \R$, which will be vital to conclude our contradiction argument.\\[6pt]
\noindent
Now $\xj \in \lj$ and $x'=0$ when $x=\xj$. So
$$ \mbox{dist} (0, \partial \Omega_j) = \kj^{\half} \mbox{dist} (\xj, \partial \Omega) \geq \beta \kj^{\xi}
\rightarrow \infty ~~~\mbox{as}~~~ \kj \rightarrow \infty.$$
Thus given an arbitrary compact subset $K$ of $\RN$, $K \subset \oj$ for $j$ sufficiently large, and hence given
$T>0$, $K \times [-T, T] \subset \oj \times [-\kj \tj, \infty)$ and is uniformly bounded away from
$\partial \oj \times \{ - \kj \tj \}$ for $j$ sufficiently large.
And since $0 \leq U^{\kj}, V^{\kj} \leq M$ for all $j$, it follows from the
interior estimates of \cite[p 342]{LSU} that 
$U^{\kj}, V^{\kj}$ are 
bounded independently of $j$ 
in $W^{2,1}_p (K \times [-T, T])$ for every $p \in [1, \infty)$ and thus in 
$C^{1+\lambda, \frac{1+\lambda}{2}}(K \times [-T, T])$ for every $\lambda \in (0,1)$
(see \cite[p 5]{LSU} for
the definition of the parabolic space $W^{2,1}_p (K \times [-T, T]))$.
So there is a subsequence of $U^{\kj}, V^{\kj}$ that
converges strongly in
$C^{1+\lambda,  \frac{1+\lambda}{2}}(K \times [-T, T])$ for each $\lambda \in (0,1)$. 
Thus since $\kj^{-1}f(U^{\kj}), \kj^{-1}g(V^{\kj}) \ra 0$ uniformly
(on $K \times
[-T,T]$) as $j \ra \infty$, passing to the limit in the weak form of
(\ref{nveq})  yields a weak solution $U,V \in C^{1+\lambda,  \frac{1+\lambda}{2}}(K \times
[-T, T])$ of the system
\begin{equation}
\label{limitsys}
\begin{array}{l}
U_t = \Delta U - UV \\
V_t = \Delta V - \alpha UV.
\end{array}
\end{equation}
That in fact $U,V \in C^{2+\lambda,  1+\frac{\lambda}{2}}(K \times[-T, T])$ and is a classical
solution of (\ref{limitsys}) on
$\mbox{int} (K \times [-T, T])$ then follows
immediately from \cite[p 224]{LSU}. And thus by a diagonalisation argument, a
subsequence of $ U^{\kj}, V^{\kj}$ converges uniformly on compact subsets of
$\RN \times \R$ to a solution $U,V$ of (\ref{limitsys}) with $0 \leq U,V \leq M$ and
$U(0,0)  \geq  \epsilon_0,  V(0,0) \geq  \epsilon_0$.\\[15pt]
\noindent
Now consider part $(ii)$,  for which condition \bfive~ is assumed to hold. 
Proceeding by contradiction as for part $(i)$,  the argument above leading to the system (\ref{nveq})
follows through with the single change that now $\xj \in \Omega$ instead of $\xj \in \lj$.
This leads to  there being two possible types of limit problem that arise from letting $j \rightarrow
\infty$ in (\ref{nveq}) in this case, depending on the behaviour of the
sequence $\{ \xj \}_{j=1}^{\infty}$. If $\mbox{dist}(0, \partial \oj) \rightarrow \infty$
for a subsequence as $j \rightarrow \infty$, then exactly as above, we obtain a solution
$(U,V)$ of (\ref{limitsys}) on $\RN \times \R$.
The second possible type of limit problem arises if \{dist$(0,
\partial \Omega_j)\}_{j=1}^{\infty}$ is bounded. In this case, there is a
subsequence (not re-labelled) for which $\Omega_j$ approaches a
half-space $H$ as $j \rightarrow \infty$ in the sense of the
definition below. There exists a subsequence of $\{x_j\}_{j=1}^{\infty}$
which we denote again by $\{x_j\}_{j=1}^{\infty}$ and a point $x_0$ such that
$$
x_j \rightarrow x_0 \qquad \mbox{as} \qquad j \rightarrow \infty.
$$
Since by the rescaling,
\begin{equation}
\label{dist} \mbox{dist} (0, \partial \Omega_j) = k_j^{1/2} \mbox{dist} (x_j,
\partial \Omega),
\end{equation}
and since by hypothesis $\mbox{dist}(0, \partial \Omega_j)$ is
bounded independently of $j$, (\ref{dist}) implies that $x_0 \in
\partial \Omega$. Furthermore, it turns out that $\partial H$ is parallel
to the tangent plane to $\partial \Omega$ at $x_0$. 
The precise sense of the convergence of the sequence
$\{\Omega_j\}_{j=1}^{\infty}$ is as follows.\\

\noindent {\bf Definition.} We say that $\Omega_j$ approaches a
half-space $H$ as $j \rightarrow \infty$ if : (i) let $K \subset
H$ be an arbitrary compact set contained in H; then there exists
$\jz$, depending on $K$, such that $K \subset \Omega_j$ for all $j
\geq \jz$; (ii) similarly let $K' \subset ~\mbox{int}(\RN
\setminus H)$ be an arbitrary compact set; then there exists
$\jz'$ depending on $K'$ such that $K' \subset \RN \setminus
\Omega_j$ for all $j \geq \jz'$.\\

\noindent Note that since $0 \in \Omega_j$ for each $j$, $0 \in
\overline{H}$ and it follows as above that as $\{k_j\} \rightarrow
\infty$, $(U^{\kj}, V^{\kj})$ converges uniformly on compact
subsets of $H \times \R$ to a function pair $(U,V)$, which is
continuous on $H \times \R$ and solves (\ref{limitsys}) on $H
\times \R$. The limiting pair $(U,V)$ is in fact {\em uniformly}
continuous on $\overline{H} \times [-T,T]$ for each $T>0$.
This follows from the fact that for some $\lambda >0$, $(U^{k_j}, V^{k_j})$
is bounded in $C^{\lambda, \frac{\lambda}{2}}(\overline{\Omega_j} \times [-T, T])$
independently of $j$ for each $T>0$, which can be proved in the following steps:
\begin{itemize}
\item[(i)] straightening the boundary of $\Omega_j$ locally (as done, for instance, 
in \cite[p97-98]{GT}) leads to transforming (\ref{nveq}) into a more complicated system for
a transformed pair $(\tilde{U}^{k_j}, \tilde{V}^{k_j})$ which can be shown to be bounded in the
$C^{\lambda, \frac{\lambda}{2}}$-norm of its flat domain independently of $j$ for some
$\lambda >0$ using \cite[p204]{LSU};
\item[(ii)] reversing the straightening of the boundary then gives the required uniform H\"older
bound on $(U^{k_j}, V^{k_j})$ because the $C^{2,\mu}$ norm of the function which
defines the boundary at a point of $\partial \Omega_j$ is bounded
from above independently of $j$ by the $C^{2,\mu}$ norm of $\partial \Omega$.
\end{itemize}

%

\noindent We now claim that $U=m_1(x_0)$ and $V=m_2(x_0)$ on $\partial H \times \R$,
where $m_1, m_2$ are as in \bthree. To see this, first note that
by Lemma \ref{lemA}, $ m_i^{k_j}(x_j) \rightarrow m_i(x_0) \geq 0$ for $i
\in \{1,2 \}$. Now let $y \in \partial H$. Since $\Omega_j$
converges to $H$ as $j \rightarrow \infty$, there is a sequence
$\{ s_j \}_{j=1}^{\infty}$ with $s_j \in \partial \Omega_j$ such
that $s_j \rightarrow y$ as $j \rightarrow \infty$. And
$U^{k_j}(s_j, t) = m_1^{k_j}(x_j+\frac{s_j}{\sqrt{k_j}})
\rightarrow m_1(x_0)$ as $j \rightarrow \infty$, by Lemma \ref{lemA}. To
see that $U(y,t)=m_1(x_0)$, fix $\tilde{y} \in \mbox{int}H$ and let $j_1$
be such that for $j \geq j_1$, $\tilde{y} \in \Omega_j$ and
$\|\tilde{y}-s_j \| \leq 2\| \tilde{y}-y\|$. Then since for each
$t \in \R$, $\{ U^{k_j}(\cdot, t)\}_{j=1}^{\infty}$ is
equicontinuous on $\overline{\Omega}_j$, given $\epsilon >0$,
there exists $\delta >0$, independent of $j$, such that
$$|U^{k_j}(\tilde{y},t)-\mkjone(\xj + \frac{s_j}{\sqrt{\kj}})|< \epsilon ~~\mbox{if}~~j \geq j_1 ~\mbox{and}~
2\|\tilde{y}-y\|< \delta,$$ and letting $j \rightarrow \infty$
gives that $$|U(\tilde{y}, t) - m_1(x_0)| \leq \epsilon
~~~\mbox{if}~~~ 2 \| \tilde{y} - y \| < \delta.$$ Letting
$\tilde{y} \rightarrow y$, it follows that $U(y,t)=m_1(x_0)$, and
similarly, that $V(y,t)=m_2(x_0)$, as required.\\[6pt]
\noindent 
So we have a solution $U,V$ of (\ref{limitsys}) on $H \times \R$ with
$0 \leq U,V \leq M$. Moreover, by condition \bfive, at least one of $U, V$ is identically zero
on $\partial H$. And as in the first possible limit problem above,
$U(0,0) \geq \epsilon_0$
and $V(0,0) \geq \epsilon_0$.\\[15pt]
\noindent 
Thus a contradiction approach to proving both parts $(i)$ and $(ii)$
of Lemma \ref{keylem} leads to a solution
of the limit equations (\ref{limitsys}) on either $\RN \times \R$ or on $H \times \R$ for a
half-space $H \subset \RN$.
In what follows we complete the proof by showing that for both possible limit problems, at
least one of $U,V$ must be identically zero, which is inconsistent with $U(0,0) \geq \epsilon_0$
and $V(0,0) \geq \epsilon_0$. We focus on the details of the case where $U,V$ are
defined on $H \times \R$; the case when $U,V$ are defined on $\RN \times \R$ is
slightly simpler and is treated in \cite{DZ}.\\[6pt]
\noindent
Note first that $U$ and $V$ are constant on $\partial H \times \R$.
We can suppose that $H= \{ x : x_N >0 \}$ without loss of generality, since $\Delta$
is invariant under rotation and translation of the spatial domain.
Now extend $\eta := \alpha U-V - (\alpha U|_{\partial H}-V|_{\partial H})$
to a function $\hat{\eta}$ on $\RN \times \R$ which is odd about
$\partial H$ in the direction orthogonal to $\partial H$, so that for
$(x,t)$ with $x_N <0$,
\begin{equation}
\label{extend}
\eh(x_1, \ldots, x_{N-1}, x_N, t) =
- \eta (x_1,\ldots, x_{N-1}, - x_N, t).
\end{equation}
It follows immediately from (\ref{limitsys}) that on
$\{x_N>0\} \times \R$, $\eh$ is
pointwise classically differentiable up to second order in space and first
order in time and $\eh_t = \Delta \eh$.
And the extension construction (\ref{extend}) gives that the same holds
in $\{ x_N <0 \} \times \R$, since $\eh_t = \Delta \eh$
is autonomous and all spatial derivatives are of even order.
Now let $\phi \in C_0^{\infty} (\RN \times \R)$ be supported in a ball
${\mathcal{B}}$ in $\RN \times \R$, and note that the outward
unit normals $\nu, \tilde{\nu}$ to $\{x_N>0\} \times \R$ and
$\{x_N<0\} \times \R$ respectively are the $(N+1)$-vectors
$\nu = (0, \ldots, 0,-1,0)$ and $\tilde{\nu} = (0, \ldots, 0,1,0)$.\\[6pt]
\noindent
 Then for each $i = 1, \ldots, N$, Green's Theorem gives that

\begin{eqnarray*}
\lefteqn{ \int_{{\scriptstyle \RN \times \R}} \eh \phi_{x_i} = \int_{{\scriptscriptstyle {\mathcal{B}} \cap (\{x_N>0\} \times \R)}} \eh \phi_{x_i} +  \int_{{\scriptscriptstyle {\mathcal{B}} \cap (\{x_N<0\} \times \R)}} \eh \phi_{x_i}} \\
 & = & - \int_{{\scriptscriptstyle {\mathcal{B}} \cap (\{x_N>0\} \times \R)}} \eh_{x_i} \phi + \int_{{\scriptscriptstyle {\mathcal{B}}  \cap (\{x_N=0\} \times \R)}} \eh \phi \nu_i - \int_{{\scriptscriptstyle {\mathcal{B}} \cap (\{x_N<0\} \times \R)}} \eh_{x_i} \phi + \int_{{\scriptscriptstyle {\mathcal{B}} \cap (\{x_N=0\} \times \R)}} \eh \phi \tilde{\nu}_i\\
 & = & - \int_{{\scriptstyle \RN \times \R}} \eh_{x_i} \phi
\end{eqnarray*}
since $\nu = - \tilde{\nu}$ and so the boundary terms cancel.
Thus $\eh$ has weak first order spatial (and likewise, first order time and
second order spatial) derivatives that equal
the pointwise classical derivatives away from $\partial H \times \R$.
So $\eh$
is a weak solution of $\eh_t = \Delta \eh$ on any bounded subdomain of
$\RN \times \R$.
And since $\eta$ is
continuous on $\overline{H} \times \R$ (since $U$ and $V$ are)  we have that
$\hat{\eta}$ is continuous on $\RN \times \R$, and hence by \cite[p223, Thm 12.1]{LSU},
$\hat{\eta}$ is a classical solution of $\hat{\eta}_t = \Delta \hat{\eta}$
on $\RN \times \R$.
So the fact that a bounded
solution of $\hat{\eta}_t = \Delta \hat{\eta}$ on $\RN \times \R$
must be constant \cite{John} implies that $\hat{\eta} \equiv
\hat{\eta}|_{\partial H} =0$. Thus
 on $\overline{H} \times \R$,
either
$\alpha U - V = \alpha m_1(x_0)$, if $V=0$ on $\partial H$, or
$\alpha U - V = -m_2(x_0)$, if $U=0$ on $\partial H$.\\[6pt]
\noindent
Consider the case when $\alpha U - V = -m_2(x_0)$  on $\overline{H} \times
\R$ (a similar argument
applies if $\alpha U-V= \alpha m_1(x_0)$). Then  $U$
satisfies
\begin{equation}
\label{sub}
\begin{array}{ll}
U_t=\Delta U - U(\alpha U + m_2(x_0)) & \mbox{on}~H \times \R, \\
U=0 &  \mbox{on}~\partial H \times \R.
\end{array}
\end{equation}
We will show that $U \equiv 0$.
Note that, unlike in the homogeneous Neumann boundary condition case considered
in \cite{DZ}, here it is necessary to consider (\ref{sub}) on $H \times \R$
rather than on $\RN \times \R$ because our extended function $\hat{\eta}$ is odd rather
than even.\\[6pt]
\noindent
Since $0  \leq U \leq M$ and $U$ is constant on $\partial H \times \R$, well-known local estimates
\cite{LSU} imply that $U$ is bounded in $C^{2+ \lambda,1 + \frac{\lambda}{2}}$ uniformly in
$\overline{H} \times \R$ for each $\lambda \in (0,1)$.
Define $z(t)=\sup_{x \in \overline{H}}U(x,t)$. To see that $z$ is Lipschitz
and thus differentiable almost everywhere (\cite[p 81]{EG}),
take $s,t \in \R$ and let $a_n \in H$ be such that $U(a_n, t) \geq z(t) - 1/n$.
Then
$$ z(t) - z(s) \leq U(a_n, t) + 1/n - \sup_{x \in \overline{H}} U(x,s)
               ~  \leq ~  U(a_n, t) - U(a_n, s) + 1/n
              ~ \leq ~ M |t-s| + 1/n $$
for some $M>0$, since $U$ is bounded in $C^{2,1}$ uniformly in $\bar{H} \times \R$.
Similarly, $z(s) - z(t) \leq M|t-s|$. \\[6pt]
\noindent
Now fix $\tb \in \R$. Because $U \geq 0$ on $\overline{H}\times \R$, $U=0$
on $\partial H \times \R$ and $U$ is {\em uniformly} continuous on $\overline{H} \times \{ \tb \}$, $U(\cdot,
\tb)$ either attains its supremum over $x \in \overline{H}$ at some $\xb \in
\mbox{int}H$ or there exists a sequence $\xn$ with dist$(\xn, \partial H) \geq \delta >0$
for every $n$ and $U(\xn, \tb) \to \sup_{x \in \overline{H}}U(x,t)$ as $n
\to \infty $. Suppose that $U(\cdot, \tb)$ attains its supremum at $\xb \in
\mbox{int}H$. It follows from the definition of $z$ that for $h>0$,
$$h^{-1}(z(\tb) - z(\tb - h)) \leq h^{-1}(U(\xb, \tb) - U(\xb, \tb - h )).$$
Hence, since $\Delta U(\xb, \tb) \leq 0$ and $m_2(x_0) \geq 0$,
\begin{eqnarray*}
\mbox{limsup}_{h \to 0^+} \frac{z(\tb ) - z(\tb-h)}{h} & \leq & U_t(\xb, \tb) \\
& \leq & -U(\xb,\tb) (\alpha U(\xb, \tb) + m_2(x_0)) \\
& \leq & -\alpha z(\tb)^2.
\end{eqnarray*}
If  $\sup_{x \in \overline{H}}U(x,\tb)$ is not attained, let $\xn$ be such that
$U(\xn, \tb) \to \sup_{x \in \overline{H}}U(x,\tb)$ as $n\to \infty $.
Now the local estimates \cite{LSU} clearly imply that
a subsequence
of $U(\cdot + \xn, \cdot)$ converges uniformly on compact sets of either $\RN \times
\R$ if dist$(\xn, \partial H) \to \infty$, or else $H \times \R$ for
some (possibly different) half-space
$H$, to a solution $\ut$ of (\ref{sub}).
Note that in both cases, $0$ belongs to the interior of the
domain of $\ut$,  that
$$ \ut(0, \tb) = \lim_{n \to \infty} U( \xn, \tb)= \sup_x U(x, \tb) =: z(\tb)$$
and also that
$$\ut(0, \tb) = \sup_x \ut(x, \tb).$$
\noindent Next define $\zt (t)= \sup_x \ut(x,t), t \in \R$. Then since $\ut$ satisfies
(\ref{sub}) and $\ut(\cdot, \tb)$ attains its supremum in the interior of the
domain of $\ut$, the argument given above in the analysis of $z$ applies to
$\zt$ to give
$$ \mbox{limsup}_{h \to 0^+} \frac{ \zt ( \tb) - \zt(\tb-h)}{h}
\leq - \alpha \zt (\tb)^2.$$
\noindent
Now, we have immediately that $\zt(\tb)= \ut(0, \tb)= z(\tb)$. And for $t \in \R$,
$\zt(t) \leq z(t)$,
 since if $\zt(t)> z(t)$ for some $t$, then $\sup_x \ut(x,t) >
\sup_x U(x,t)$, so there exists $\xtt$ with $\ut (\xtt, t) > \sup_x U(x,t)$,
and since $U(\xtt + \xn, t) \to \ut(\xtt, t)$, there exists $n_0$ for which
$U(\xtt + x_{n_0}, t) > \sup_x U(x,t)$, which is impossible. 
So
$$ \mbox{limsup}_{h \to 0^+} \frac{ z( \tb ) - z(\tb-h)}{h}
\leq \mbox{limsup}_{h \to 0^+} \frac{ \zt( \tb ) - \zt(\tb-h)}{h}
\leq - \alpha \zt (\tb)^2 \leq - \alpha z(\tb)^2.$$
It follows that for every $t \in \R$, $\mbox{limsup}_{h \to 0^+}
{\displaystyle \frac{z( \tb) - z(\tb-h)}{h}}
\leq - \alpha z(\tb)^2$, and hence on the set of full measure on which $z$ is
differentiable,
\begin{equation}
\label{diffeq}
\dot{z}\leq - \alpha z(t)^2.
\end{equation}
If there exists $t$ with $z(t) >0$, then $z(s) \geq z(t)>0$ for all $s \leq t$, since
(\ref{diffeq})
implies that $z$ is non-increasing. So since $z(s) \leq M$, $\dot{z} /z(t)^2
\leq \dot{z} / z(s)^2 \leq \dot{z} /M^2$ for $s \leq t$. Thus $\dot{z} /z^2 \in L^1(t_0, t)$, $t_0
< t$. So
(\ref{diffeq}) can be integrated to obtain that for any $t_0 < t \in \R$,
$1 /z(t) \geq \alpha (t-t_0) + 1/z(t_0)$
and so
\begin{equation}
\label{est}
z(t) \leq \alpha^{-1}(t-t_0)^{-1}.
\end{equation}
Since $t_0 \in \R$ was arbitrary, we can
let $t_0 \to - \infty$ in (\ref{est}) to find that $z(t)=0$ for every $t \in
\R$. Hence $U \equiv 0$.\\[6pt]
\noindent
If $\alpha U - V = \alpha m_1(x_0)$, a similar argument involving substitution for
$U$ in the $V$ equation implies that $V \equiv 0$.
\qed

\noindent Recall the definition of $\wk$ from (\ref{wk}).
Lemma \ref{keylem} yields the following convergence result, which gives
 the convergence properties
that will be used in subsequent sections to
analyse the long-time behaviour of solutions of problem $\pk$ for large $k$.

\begin{lemma}
\label{convlem}
$(i)$ If $(\uk, \vk)$ is a solution of $\pk$ with $0 \leq \uk, \vk \leq M$, then for
each $\tz >0$,
\begin{equation}
\sup_{t \geq \tz} \left\{ \| (\wk)^+ - \alpha \uk \|_{L^2(\Omega)} + \| (\wk)^- + \vk \|_{L^2(\Omega)} \right\}
\rightarrow 0   ~~~\mbox{as}~~k \rightarrow \infty
\label{conv1}
\end{equation}
and
\begin{equation}
\label{wkeq}
\wk_t = \Delta \wk + h(\wk) + R(\uk, \vk)
\end{equation}
where
\begin{equation}
\label{defh}
h(w) := \alpha f(\alpha^{-1}w^+) - g(-w^-),
\end{equation}
and
\begin{equation}
\label{convR}
\sup_{t \geq \tz} \| R(\uk, \vk) \|_{L^2(\Omega)} \rightarrow 0 ~~~\mbox{as}~~k \rightarrow \infty.
\end{equation}
$(ii)$ If, in addition, $\mkone$, $\mktwo$ satisfy the supplementary condition \bfive, then (\ref{conv1})
and (\ref{convR}) hold with the norm $\| \cdot \|_{L^2(\Omega)}$ replaced by the norm
$\| \cdot \|_{L^{\infty}(\Omega)}$. \\ \\
(Here $w^+ := \max \{ w,0 \}$, $w^- := \min \{ w,0 \}$ and thus $w= w^+ +  w^-$.)
\end{lemma}

\pf Fix $\tz >0$ and let $\epsilon >0$. Lemma \ref{keylem} implies that there exists $k_0$
such that for each $k\ \geq k_0$, $x \in \lk$ and $t \geq \tz$, either
$\alpha \uk (x,t) \leq \epsilon$ or $\vk (x,t) \leq \epsilon$.\\[6pt]
\noindent
Suppose first that $\alpha \uk (x,t) \geq \ep$ and $\vk \leq \ep$. Then since $\wk =
\alpha \uk - \vk \geq 0$, $\wkp = \wk$ and $\wkm = 0$.
So $| \wkp - \alpha \uk | = | \vk | \leq \ep$ and $| \wkm + \vk| = |\vk| \leq \ep$.
Similarly, if $\alpha \uk (x,t)\leq \ep$ and $\vk (x,t)
 \geq \ep$, then $\wk \leq 0$, $\wkp = 0$ and $\wkm = \wk$. And hence
$| \wkp - \alpha \uk | = |\alpha \uk| \leq \ep$ and $| \wkm + \vk| = |\alpha \uk| \leq \ep$.
Finally, if $\alpha \uk (x,t)\leq \ep$ and $\vk (x,t) \leq \ep$, then
$| \wkp - \alpha\uk|
 \leq |\wk| + |\alpha \uk| \leq 3 \ep$ and $| \wkm + \vk| \leq 3 \ep$.\\[6pt]
\noindent
Lemma \ref{keylem} ensures that one of these three possibilities must arise for each
$(x,t) \in \lk \times [\tz, \infty)$ where $k \geq k_0$. So for such $(x,t)$ and $k$,
\begin{equation}
\label{conv1ep}
 | \wkp (x,t) - \alpha\uk(x,t)| \leq 3 \ep ~\mbox{and}~ | \wkm (x,t) + \vk(x,t)| \leq
3\ep,
\end{equation}
and hence for $R: \R \times \R \rightarrow \R$ defined by
$$ R(\uk, \vk) := \alpha f(\uk) - g(\vk) - \alpha f(\alpha^{-1} (\wk)^+)+g(-(\wk)^-),$$
we have
\begin{eqnarray}
|R(\uk, \vk)| & \leq &  \alpha K_f | \uk -
\alpha^{-1}\wkp | + K_g |\vk + \wkm| \nonumber \\
& \leq & 3( K_f + K_g) \ep, \label{convRep}
\end{eqnarray}
where $K_f, ~K_g$ are the Lipschitz constants of $f,g$ respectively restricted
to the interval $[-M-1, M+1]$.
And since $0 \leq \uk, \vk \leq M$ and the measure $| \Omega \setminus \lk | \rightarrow
0$ as $k \rightarrow \infty$, it follows that
$$ \int_{\Omega \setminus \lk} | (\wk)^+ - \alpha \uk |^2 + | (\wk)^- + \vk |^2 \rightarrow
0~~~\mbox{and}~~~ \int_{\Omega \setminus \lk} |R(\uk, \vk)|^2 \rightarrow 0$$
as $k \rightarrow \infty$ uniformly in $t \geq \tz$. Since $| \Omega |< \infty$,
this, together with (\ref{conv1ep}) and (\ref{convRep}) establishes (\ref{conv1}) and (\ref{convR}).\\[6pt]
\noindent
If, in addition, condition \bfive~ holds, then it follows from Lemma \ref{keylem} that
the above argument holds with $\lk$ replaced by $\Omega$ throughout, from which the
last statement of Lemma \ref{convlem} is immediate. \qed

\noindent {\bf Remark} We conclude this opening section by noting a relation with
a special case of \cite{CDHMN}. In \cite{CDHMN}, a spatial segregation limit is
derived for the generalisation of problem $\pk$
in which the diffusion coefficients of $u$ and $v$ are allowed to differ. It is shown that
$\uk \rightarrow u$ and $\vk \rightarrow v$ in $L^2( \Omega \times (0, T))$ for every $T>0$,
where $uv =0$ almost everywhere in $\Omega \times (0, T)$ and $w= \alpha u -v$ is the
unique weak solution of a limiting free boundary problem (see \cite[Section 3]{CDHMN} for details).\\[6pt]
\noindent
Now if the diffusion coefficients are in fact the same, then
for each $0 < \tz < T$,
\begin{equation}
\label{convCDHMN}
\wk= \alpha \uk - \vk  \rightarrow w ~~\mbox{in}~~ C^{1+\lambda^{\prime}, \frac{1+ \lambda^{\prime}}{2}}
(\overline{\Omega} \times [\tz, T]) ~~\mbox{for all}~~\lambda^{\prime} \in (0, \lambda).
\end{equation}

\section{Long-time behaviour (1) : closeness to stationary solutions of the limit problem}
In this section we will show that for sufficiently large $k$, solutions of $\pk$ are close
to stationary solutions of a certain limit problem for sufficiently large time.
%
The appropriate notion of limit-problem stationary solutions is as follows.
Recall the definition of $h$ from Lemma \ref{convlem} and note from \bthree~ that $m_1, m_2 \in \wtp$.
We will say that $w \in \wtp $
is a solution of $\stat$ if
\begin{equation}
\label{stationary}
(S) \left\{ \begin{array}{ll} \Delta w + h(w) =0 & ~~\mbox{in}~~\Omega,\\
w = \alpha m_1 - m_2 & \mbox{on}~~~ \partial \Omega,
\end{array} \right.
 \end{equation}
 which immediately implies that $w \in C^2(K)$ for all compact sets $K \subset \Omega$
 and that $w \in C^{1, \lambda}(\overline{\Omega})$. 
Note that the results in this and the following section hold whether or not the supplementary condition
\bfive~ holds.
\bigskip

We  first collect some standard regularity, boundedness and compactness results for solutions
$\ukvk$ of $\pk$ and $\wk$ that will be useful in this and the  following sections.

\begin{lemma}
\label{lemregandcompact}
Let $\ukvk$ be the solution of $\pk$ for some $k \in \N$.
\begin{description}
\item[(a)] $u_{xt}^k (x,t)$, $v_{xt}^k (x,t)$  exist for each $\xt \in \Omega \times (0, \infty)$,
and there exists $\lambda >0$ such that for each $t>0$, $u^k_t \fixt, v^k_t \fixt \in \conel$;
\item[(b)] for each $\tz >0$, there exists $M_k$ (dependent on $k$) such that for each $t \geq t_0$,
$$ \| \uk \fixt \|_{\wtp},  \| \vk \fixt \|_{\wtp} \leq M_k, $$
and there exists $\tilde{M}$ (independent of $k$) such that for each $t \geq \tz$,
$$ \| \wk \fixt \|_{\wtp} \leq \tilde{M}, $$
(note that since $p>N$, these estimates clearly also hold with $\wtp$ replaced by $C^{1,\lambda}(\Omega)$ for some
$\lambda >0$);
\item[(c)] given $\tz >0$, there exist compact subsets
$\tilde{\Lambda}_{\tz, k} \subset \wtp \times \wtp$ (dependent on $k$) and
$\Lambda_{\tz} \subset \wtp$ (independent of $k$) such that for all $t \geq \tz$,
$$ \ukvk \fixt \in  \tilde{\Lambda}_{\tz, k} ~~~\mbox{and}~~~ \wk \fixt \in \Lambda_{\tz};$$
\item[(d)] the $\omega-$limit set (omega-limit set) $\Gamma$ of
$(u^k_0, v^k_0)$ in $\wtp \times \wtp$ is non-empty, compact, invariant, connected. and
$dist((\uk, \vk) \fixt, \Gamma) \rightarrow 0$ as $t \rightarrow \infty$, where $dist$
is measured in the $\wtp \times \wtp$ norm and, as usual,
\begin{eqnarray*} \lefteqn{ \Gamma  =  \{ (u,v) \in \wtp \times \wtp :} \\
& & ~~\mbox{there exist}~~ t_n \rightarrow \infty~~
\mbox{such that}~~ \| \uk (\cdot, t_n) - u \|_{\wtp} +  \| \vk (\cdot, t_n) - v \|_{\wtp}
\rightarrow 0 \}.\end{eqnarray*}
\end{description}
\end{lemma}

\pf We use the semiflow framework of \cite[Chp 3]{Henry} in the space $X = \lp \times \lp$ with the
domain $\wtpz \times \wtpz$ for the system with homogeneous boundary conditions discussed in the
proof of Lemma \ref{global} (respectively $X= \lp$, domain $ \wtpz$, for the corresponding homogeneous equation for the linear combination
$\wk$).

Part {\bf (a)} follows from \cite[Thm 3.5.2]{Henry} and the fact that given
$\beta <1$ sufficiently close to $1$, the fractional power space $X^{\beta}
\subset \conel$ for some $\lambda >0$.
That for such a $\lambda>0$ $\| u^k_t \fixt \|_{\conel}$, $\| v^k_t \fixt \|_{\conel}$ are
bounded independently of $t \geq \tz$ for each fixed $k$, and $\| \wk_t \fixt \|_{\conel}$
is bounded independently of $t \geq \tz$ and $k \in \N$, follow from \cite[Thm 3.5.2]{Henry}
and Lemma \ref{apriori}
(note that Lemma \ref{apriori} and the remark following (\ref{wkeqn}) give the independence
of $k$ of the bound on $\| \wk_t \fixt \|_{\conel})$.
Part {\bf (b)} then follows using Lemma \ref{apriori} again, together
with \cite[Lem 9.17]{GT}.
For Part {\bf (c)}, note from the bounds on $\uk_t, \vk_t, \wk_t$ just observed, together
with Part {\bf (b)} and condition {\bf (a)} on $f$ and $g$, that given $\tz >0$, there
are compact subsets $\tilde{\Lambda}^{\prime}_{\tz, k} \subset L^p (\Omega) \times L^p (\Omega)$
(dependent on $k$) and $\Lambda^{\prime}_{\tz} \subset L^p (\Omega)$ (independent of $k$)
such that for all $ t \geq \tz$,
$ \uk_t \fixt$, $ \vk_t \fixt$, $(f(\uk) - k \uk \vk ) \fixt$, $(g(\vk) - \alpha k \uk \vk ) \fixt$
$\in  \tilde{\Lambda}^{\prime}_{\tz, k}$ and
$ \wk_t \fixt \in \Lambda^{\prime}_{\tz};$ that $ \Lambda^{\prime}_{\tz}$ can be chosen so that
we also have $(\alpha f(\uk) - g(\vk)) \fixt \in \Lambda^{\prime}_{\tz}$ follows using
(\ref{wkeq}), (\ref{defh}), (\ref{convR}) from Lemma \ref{convlem}, in addition to Part {\bf (b)}
and condition {\bf (a)}.
Thus $\Delta \uk \fixt$, $\Delta \vk \fixt \in \tilde{\Lambda}^{\prime}_{\tz, k}$
and $\Delta \wk_t \fixt \in \Lambda^{\prime}_{\tz}$ and Part {\bf (c)} follows  since
$\Delta$ maps $\wtp$ bijectively to $L^p (\Omega)$.
Part {\bf (d)} is then immediate from {\bf (c)} and \cite[Thm 4.3.3]{Henry}. \qed

\medskip
\noindent The main result of this section is the following.
\medskip

\begin{theorem}
\label{energythm}
Suppose that solutions of $\stat$ are isolated in $L^2 (\Omega)$. Then given $\epsilon >0$ and $M>0$,
there exists $k_0$ such that for $ k \geq k_0$ and $(\uk, \vk)$ the solution of $\pk$ with $0 \leq \uk, \vk \leq M$,
there exists a solution $\tw$ of $\stat$ such that
\begin{equation}
\label{closeenergy}
\| \alpha \uk (\cdot,t) - \tw^+(\cdot)\|_{L^2(\Omega)} + \| \vk (\cdot,t) + \tw^-(\cdot)\|_{L^2(\Omega)}
\leq \epsilon
\end{equation}
for all $t$  sufficiently large (where how large $t$ needs to be depends on
$k$).
\end{theorem}

\pf First some preliminary remarks.\\[1pt]

\noindent {\em (i)} We will show that there exists a solution $\tw$ of $\stat$ such that
$\wk = \alpha \uk - \vk$ is close to $\tw$ in $L^{2}(\Omega)$ for large time.
(\ref{closeenergy}) will follow from this together with (\ref{conv1}) from Lemma \ref{convlem}.\\[1pt]

\noindent {\em(ii)} Denote the set of all solutions of $\stat$ by $\ms$. Then $\ms$ is a compact
subset of $L^2(\Omega)$ since $\ms$ is bounded in $W^{1,2}(\Omega)$ and thus any sequence in
$\ms$ has a subsequence that converges weakly in $W^{1,2}(\Omega)$ and $L^2(\partial \Omega)$,
and strongly in $L^2(\Omega)$, to a limit $w$ which is a solution
of the weak form of $\stat$, and hence, by regularity, of $\stat$.
Hence $\ms$ is a {\em finite} set, by the assumption that the solutions of $\stat$ are isolated in
$L^2(\Omega)$. In the rest of the proof, let $$\ms = \{ \wsi : 1 \leq i \leq r, ~ r \in \N \}.$$

\noindent {\em (iii)} It follows from Lemma \ref{lemregandcompact} (c)
that for each $\tz >0$,
$\wk (\cdot, t)$ lies in compact subsets of $W^{1,2}(\Omega)$ and
$C(\overline{\Omega})$ independently of $k$ and of $t \geq \tz$.
Since a continuous bijection on a compact set is a homeomorphism, 
this implies that the $L^{\infty}$-, $L^2$- and $W^{1,2}$-norms generate
equivalent metrics on the set $\{ \wk (\cdot, t) : k \in \N ~\mbox{and}~t \geq t_0 \}$.\\[2pt]



\noindent {\em (iv)} Fix $\eta >0$ and $t_0 >0$. Then there exist $\delta >0$ and $k_0$ such that
if for $k \geq k_0$ and $t \geq t_0$,
\begin{equation}
\label{boundbelow}
\begin{array}{l} \| \wk (\cdot, t) - \wsi \|_{L^2(\Omega)} \geq \eta \\
~~\mbox{for every}~~ 1 \leq i \leq r \end{array} ~~~\Rightarrow ~~~
\| \Delta \wk (\cdot, t) + h(\wk (\cdot, t)) \|_{L^2(\Omega)} \geq \delta.
\end{equation}
For if not, there exist  sequences $k_j$ and $t_j \geq t_0$
such that $k_j \rightarrow \infty$ as  $j \rightarrow \infty$, $w_{k_j}(\cdot, t_{k_j}) \in W^{2,p}(\Omega)$,
$ w_{k_j} = \alpha m_1^{k_j} - m_2^{k_j} ~~\mbox{on}~~\partial \Omega$,
$\| w_{k_j}(\cdot, t_{k_j}) - \wsi \|_{L^2(\Omega)}
\geq \eta$ for each $1 \leq i \leq r$ and $\| \Delta w_{k_j}(\cdot, t_{k_j}) +
h(w_{k_j}(\cdot, t_{k_j})) \|_{L^2(\Omega)} \rightarrow 0$ as $k_j \rightarrow \infty$. 
Then by Lemma \ref{lemregandcompact} (c), there exists $w \in W^{1,2}(\Omega)$ such
that a subsequence
$w_{k_j} \rightarrow w$ in both $W^{1,2}(\Omega)$ and $L^2(\partial \Omega)$.
Hence, using {\bf (b3)} and regularity theory, $w$
is a solution of $\stat$. But this contradicts that $\| w_{k_j}(\cdot, t_{k_j}) - \wsi \|_{L^2(\Omega)}
\geq \eta$ for each $1 \leq i \leq r$, and (\ref{boundbelow}) follows.
 \\[2pt]

Our approach, which follows closely that in \cite{DZ}, is to show that the natural energy
for the limit problem evaluated at $\wk (\cdot, t)$ decreases at a certain rate when $\wk (\cdot, t)$
lies outside $L^2$-neighbourhoods of the elements of $\ms$.
Choose and fix $\tz>0$.
For $w \in W^{1,2}(\Omega) \cap L^{\infty}(\Omega)$, define an energy
\begin{equation}
\label{energy}
\me (w) = \int_{\Omega} \frac{1}{2} | \nabla w|^2 - H(w) ~dx
\end{equation}
where $H$ is a primitive of $h$.
Note first that, by Lemma \ref{lemregandcompact} (b), $\wk (\cdot, t)$ lies in bounded sets
in $W^{1,2}(\Omega)$ and $L^{\infty}(\Omega)$ independently of $k$ and $t \geq \tz$, and hence
$\me (\wk(\cdot, t))$ is bounded independently of $k$ and of $t \geq \tz$.
Now fix $\epsilon >0$. By Lemma \ref{convlem},
\begin{equation}
\label{temp}
\wk_t = \Delta \wk + h(\wk) + R(\uk, \vk), ~~(x,t) \in \overline{\Omega} \times (0, \infty),
\end{equation}
where $\sup_{t \geq \tilde{t}_0} \| R(\uk, \vk) \|_{L^2(\Omega)} \rightarrow 0 $ as $k \rightarrow \infty$ for
each $\tilde{t}_0 >0$. And by Lemma \ref{lemregandcompact} (a), $\wk$ is sufficiently smooth that for $t \geq t_0$,
$\me (\wk (\cdot, t))$ is differentiable with respect to $t$, and
\begin{eqnarray}
\lefteqn{ \frac{d}{dt} \me (\wk (\cdot, t))  =  \int_{\Omega} \nabla \wk (x,t) \frac{\partial}{\partial t} \nabla \wk (x,t) - h(\wk (x,t)) w^k_t (x,t) ~dx \nonumber }\\
& = & \int_{\Omega} (- \Delta \wk - h(\wk)) w^k_t ~dx + \int_{\partial \Omega} \nabla \wk w^k_t ~dx
\nonumber \\
& = & \int_{\Omega} (- \Delta \wk - h(\wk)) w^k_t ~dx, \nonumber \\
& = & -  \int_{\Omega} (\Delta \wk +h(\wk))(\Delta \wk + h(\wk) +R(\uk, \vk) ), ~~~\mbox{by (\ref{temp})},
\nonumber \\
& \leq & - \| \Delta \wk + h (\wk) \|_{L^2(\Omega)} ( \| \Delta \wk + h (\wk) \|_{L^2(\Omega)}  -
\|R(\uk (\cdot, t), \vk (\cdot, t) ) \|_{L^2(\Omega)} ) \label{extra}
\end{eqnarray}
where $  \|R(\uk, \vk)\|_{L^2(\Omega} \rightarrow 0$ as $k \rightarrow \infty$ uniformly for
$t \in [t_0, \infty)$.

So if $t \geq \tz$ and
$\| \wk (\cdot, t) - \wsi \|_{L^2(\Omega)} \geq \epsilon/4$ for every $i \in \{ 1,
\ldots, r \}$, then this together with (\ref{boundbelow}) with $\eta = \epsilon/4$ gives the existence of
$\delta_1 >0$ and $k_0 \in \N$ (larger than above if necessary) such that
\begin{equation}
\label{energydecay}
\frac{d}{dt} \me (\wk (\cdot, t)) \leq - \delta_1 ~~~\mbox{for all}~~ k \geq k_0, ~t \geq t_0,
\end{equation}
since there exists $\delta >0$ such that $\| \Delta \wk + h(\wk) \|_{L^2(\Omega)} \geq \delta $ and then
 $k_0$ can be chosen so that $\|  \Delta \wk + h(\wk) \|_{L^2(\Omega)} - \|R(\uk (\cdot, t), \vk (\cdot, t) )
\|_{L^2(\Omega)} \geq \delta /2$ for $k \geq k_0$.

Denote by $\brw$ the ball in $L^2(\Omega)$, centre $w$, radius $R$. We would
like to show that there exists $\beta >0$ such that for $k$ sufficiently large,
$\me (\wk (\cdot, t))$ is a decreasing function of $t$ outside $\displaystyle{\cup_{i=1}^r
\mathcal{B}_{\beta \epsilon}(\wsi)}$ and the drop in $\me (\wk (\cdot, t))$ when $\wk (\cdot, t)$
moves from inside $\mathcal{B}_{\beta \epsilon}(\wsi)$ at some $t$ to the boundary of $\mathcal{B}_{\epsilon}
(\wsi)$ (at some later time $\tilde{t}$) is larger than the possible range of $\me (\wk (\cdot, t))$
when $\| \wk (\cdot, t) - \wsi \|_{L^2(\Omega)} \leq \beta \epsilon$. This implies that if $\wk (\cdot, t)$
moves from inside $\mathcal{B}_{\beta \epsilon}(\wsi)$ to $\partial \mathcal{B}_{\epsilon}(\wsi)$,
then $\wk (\cdot, t)$ cannot re-enter $\mathcal{B}_{\beta \epsilon}(\wsi)$ at any later time. Recall remark
{\em (iii)} and note that it follows that $\me (\wk (\cdot, t))$ is close to $\me ( \wsi)$ when $\wk (\cdot, t)$
is close to $\wsi$ in $L^2(\Omega)$ because $\me (\cdot)$ is continuous as a function of $w \in W^{1,2}
(\Omega) \cap \{ w \in L^{\infty}(\Omega) : \| w \|_{L^{\infty}(\Omega)} \leq M \}$.

To prove this, we first show that there exists $T>0$ such that if $t \geq t_0$ and $\wk (\cdot, t)
\in \mathcal{B}_{\epsilon / 4}(\wsi)$ and $\wk (\cdot, t+ \tilde{t}) \not\in \mathcal{B}_{\epsilon}(\wsi)$,
$\tilde{t} >0$, then $\tilde{t} \geq T$. To see this, recall from Lemma \ref{lemregandcompact} (b)
that $\wk (\cdot, t)$
lies in a bounded set in $W^{2,2}(\Omega)$ for all $k$ and all $t \geq t_0$, so there exists
$M_1 >0$ such that $\| \Delta \wk (\cdot, t) + h(\wk(\cdot, t)) \|_{L^2(\Omega)} \leq M_1$
for all such $k$ and $t$. And
$$ \wk (x,t_1) - \wk (x, t_2) = \int_{t_1}^{t_2} w^k_t (x,t)~dt =
\int_{t_1}^{t_2} (\Delta \wk (x,t) + h(\wk(x,t))~dt,$$
so
\begin{eqnarray*}
\int_{\Omega} | \wk (x, t_1) - \wk (x,t_2)|^2~dx & \leq & (t_2 - t_1) \int_{\Omega} \int_{t_1}^{t_2}
(\Delta \wk + h (\wk))^2~dt~dx \\
& = & (t_2 - t_1) \int_{t_1}^{t_2} \left( \int_{\Omega} ( \Delta \wk + h(\wk))^2 ~dx \right)~dt \\
& \leq & M_1^2 (t_2 - t_1)^2.
\end{eqnarray*}
Hence $t_2 - t_1 \geq \| \wk (\cdot, t_1) - \wk (\cdot, t_2) \|_{L^2(\Omega)} / M_1$, from which the
existence of $T = T_{\epsilon} = {\displaystyle \frac{3 \epsilon}{4 M_1}}$ follows.

Now this together with (\ref{energydecay}) implies that in going from $\mathcal{B}_{\epsilon / 4}(\wsi)$
to $\partial \mathcal{B}_{\epsilon}(\wsi)$, $\me (\wk (\cdot, t))$ drops by at least $T_{\epsilon} \delta_1$.
And we can choose $\beta >0 ~(\leq \frac{1}{4})$ so that for $t \geq t_0$,
\begin{equation}
\label{energyest}
\| \wk (\cdot, t) - \wsi \|_{L^2(\Omega)} < \beta \epsilon \Rightarrow | \me (\wk (\cdot, t)) - \me (\wsi)|
< \frac{1}{2} T_{\epsilon} \delta_1.
\end{equation}

Now we can apply (\ref{boundbelow}) with $\eta = \beta \epsilon$ together with (\ref{extra}) to
obtain $\tilde{k} \geq k_0$ and $\delta_2 >0$ such that for $t \geq t_0$,
$$ k \geq \tilde{k} \Rightarrow \frac{d}{dt} \me (\wk (\cdot, t)) \leq - \delta_2 $$
when $\| \wk (\cdot, t) - \wsi \|_{L^2(\Omega)} \geq \beta \epsilon$ for each $i \in \{ 1, \ldots, r \}$.
Thus for $k \geq \tilde{k}$, $\me (\wk (\cdot, t))$ decreases when $\wk (\cdot, t)$ lies
in $ \mathcal{B}_{\epsilon}(\wsi) \setminus \mathcal{B}_{\beta \epsilon}(\wsi)$, for some $i$, and the
drop is at least $T_{\epsilon} \delta_1$ as $\wk (\cdot, t)$ moves from inside $\mathcal{B}_{\beta \epsilon}(\wsi)$
to $\partial \mathcal{B}_{\epsilon}(\wsi)$, since $\beta \leq \frac{1}{4}$. It follows using
(\ref{energyest}) that if $\wk$
leaves $\mathcal{B}_{\beta \epsilon}(\wsi)$ and moves out to $\partial \mathcal{B}_{\epsilon}(\wsi)$,
it cannot re-enter $\mathcal{B}_{\beta \epsilon}(\wsi)$ at a later time.

Now if $\wk (\cdot, t) \not\in \cup_{i=1}^r \mathcal{B}_{\beta \epsilon} (\wsi)$ for all $t$
sufficiently large, then $\me (\wk (\cdot, t))$ would decrease at at least rate $- \delta_2$ for
all large time, which would contradict the fact that $\me (\wk (\cdot, t))$ is bounded below independently
of $t \geq \tz$. Hence there is a sequence of times $t_n \rightarrow \infty$ for which
$\wk (\cdot, t_n) \in \cup_{i=1}^r \mathcal{B}_{\beta \epsilon} (\wsi)$, and since there are a finite number
of $\wsi$, there exists $i_0$ and a subsequence $(t_{n_m})_{m=1}^{\infty}$ of $(t_n)_{n=1}^{\infty}$
such that $\wk (\cdot, t_{n_m}) \in \mathcal{B}_{\beta \epsilon} ( \overline{w}^{i_0})$.
But if $\wk (\cdot, t)$ left $\mathcal{B}_{\epsilon}(  \overline{w}^{i_0})$ for some $t \geq t_{n_0}$,
it would not be able to re-enter $\mathcal{B}_{\beta \epsilon}( \overline{w}^{i_0})$.
So $\wk (\cdot, t) \in \mathcal{B}_{\epsilon}(  \overline{w}^{i_0})$ for all $t$ sufficiently large.
 \qed

\section{Long-time behaviour (2) : convergence to stationary solutions of $\pk$}



\noindent Note first that here {\em all} solutions of $\stat$ are not identically equal to zero, 
since it is supposed in $\bthree$
that $\alpha m_1 - m_2$ is not identically zero on $\partial \Omega$.
Now observe that (as in \cite{DG, DZ}), it follows from \cite{CF} that such solutions
$\tw$ of $\stat$ only take the value zero on a set of measure zero. Hence $h'(\tw (x))$
exists for almost every $x \in \Omega$. This enables us to make the following definition.

\begin{definition}
\label{defnondegen}
A  solution $\tw \in \wtp$ of $\stat$ is said
to be {\em non-degenerate} if the only solution $w \in   \wtp$
of the linearised equation
\begin{equation}
\label{nondegen}
\begin{array}{ll}
\Delta w + h'(\tilde{w})w=0 & a.e.~~\mbox{in}~~ \Omega, \\
w=0 & \mbox{on}~~ \partial \Omega,
\end{array}
\end{equation}
is identically equal to zero.
\end{definition}

\begin{lemma}
\label{lemomegalimit}
Suppose that a solution $\tw$ of $\stat$ is  non-degenerate. Then given $M>0$,
there exist
$\epsilon, k_0 >0$ such that for each $k \geq \kz$, a solution $(\uk, \vk)$ of $\pk$
that is defined for all $t \in \R$ and $0 \leq \uk, \vk \leq M$ and satisfies
\begin{equation}
\label{closealltime}
\| \wk (\cdot, t) - \tw \|_{\lto} \leq \epsilon ~~~\mbox{for all}~~t \in \R,
\end{equation}
must have $\uk_t$ and $\vk_t$ identically zero on $\Omega \times \R$.\\[5pt]
(Note that for this lemma we do {\em not} assume that $\uk(\cdot, 0), \vk(\cdot, 0)$ are given
by {\bf (b4)}.)
\end{lemma}

\pf Our proof follows that in \cite[Thm 3]{DZ} which establishes the corresponding result
with zero Neumann boundary conditions. Much of the argument is un-changed
and we give an outline here,
giving most detail in a blow-up argument where the main differences with \cite{DZ} lie.

Suppose that the result is false for some $M>0$. Then there exist sequences
$\kj \rightarrow \infty$, $\epj \rightarrow 0$ and solutions $(\ukj, \vkj)$ of $(P_{\kj})$
that are defined for all $t \in \R$, $0 \leq \ukj, \vkj \leq M$,
$$\| \wkj (\cdot, t) - \tw \|_{\lto} \leq \epj~~~\mbox{for all}~~t \in \R $$
for the solution $\tw$ of $\stat$ but $( \ukj_t, \vkj_t)$ is not identically zero on
$\Omega \times \R$.

First consider $\pk$ for fixed $k$. Let $(\uk, \vk)$ be a solution defined for all $t \in \R$
with $0 \leq \uk, \vk \leq M$ and $(\uk_t, \vk_t)$ not identically zero (so $(\uk, \vk)$
is a non-stationary solution of $\pk$).

Now since $0 \leq \uk, \vk \leq M$ for all $t \in \R$, standard parabolic estimates \cite{LSU}
yield that $(\uk_t, \vk_t)$ is uniformly bounded on $\Omega \times \R$.
[Note that this bound depends, in the first instance, on $k$.]
Since $(\uk, \vk)$ is non-stationary, at least one of $\uk_t$, $\vk_t$ is non-trivial.
We introduce the norm
$$ \| (h,l) \|' = \sup_{s \in \R} (\|h \|_{L^2(\Omega \times (s, s+1))} + \| l \|_{L^2(\Omega \times (s, s+1))}) $$
which is finite for functions $h,l \in L^{\infty}(\Omega \times \R)$, in particular, for
$(\uk_t, \vk_t)$.

Now note that since $f,g$ are assumed to be continuously differentiable, bootstrapping
and differentiation gives that
$(\uk_t, \vk_t)$ is a solution of the linear system
\begin{eqnarray}
h_t & = & \Delta h + (f'(\uk) - k\vk) h - k \uk l, ~~~(x,t) \in \Omega \times \R, \label{hlsystem} \\
l_t & = & \Delta l + (g'(\vk) - \alpha k \uk) l - \alpha k \vk h, ~~~(x,t) \in \Omega \times \R \nonumber
\end{eqnarray}
with the boundary condition
$$(h,l)(x,t) =0 ~~~\mbox{for}~~~(x,t) \in \partial \Omega \times \R $$
since $(\uk, \vk)$ satisfies time-independent Dirichlet boundary conditions (by \bthree).
Since $(\uk_t, \vk_t)$ is not identically zero, we can multiply $(\uk_t, \vk_t)$ by a constant
to obtain
a solution of (\ref{hlsystem}), called $(\hk, \llk)$, say, that satisfies
\begin{equation}
\label{normalise}
\|(\hk, \llk) \|' = 1.
\end{equation}

Now a Kato-inequality argument gives that $\hk$ and $\llk$ are bounded in $L^{\infty}(\Omega \times \R)$
independently of $k$; since the proof is identical to that in \cite{DZ} {\em modulo} replacing
the zero Neumann boundary conditions for $(\hk, \llk)$ in \cite{DZ} by zero Dirichlet conditions
here, we omit the details.


We now use a blow-up argument to deduce that one of $\hkj$ and $\lkj$ is uniformly small
away from the set where $\tilde{w}=0$ if $j$ is large. More precisely, given a compact subset
$\Lambda$ of $( \mbox{int} \Omega) \setminus \{x: \tilde{w}(x) =0 \}$ and an $\epsilon_0>0$, we prove that
there exists $j_0 >0$ such that
\begin{equation}
\label{hkjlkjsmall}
| \hkj (x,t)| \leq \epsilon_0~~\mbox{or}~~| \lkj (x,t)| \leq \epsilon_0~~\mbox{if}~~(x,t) \in
\Lambda \times \R~~\mbox{and}~~j \geq j_0.
\end{equation}
Suppose that (\ref{hkjlkjsmall}) is false. Then there exist $x_j \in \Lambda$ and $t_j \in \R$
such that $|\hkj(x_j, t_j)| \geq \epz$ and $|\lkj(x_j, t_j)| \geq \epz$ for a sequence of $j's$ tending
to infinity (not re-labelled). Without loss, we can assume, by a shift in time, that $t_j =0$
for every $j$ (note that the $\ukj$ and $\vkj$ in (\ref{hlsystem}) must also be shifted in time).
Now thanks to the uniform-in-$k$ bounds on $\hkj, \lkj$ obtained above, we can rescale and blow-up
(\ref{hlsystem}) much as in the proof of Lemma \ref{keylem}. Note that since $x_j \in \Lambda$ and
$\Lambda$ is compactly contained in $(\mbox{int} \Omega) \setminus \{x: \tilde{w}(x) =0 \}$, any limit point
of the $x_j$ cannot lie on $\partial \Omega$, and hence rescaling always yields a limit system defined
on $\RN \times \R$ (rather than $H \times \R$ for a half-space $H \subset \RN$). Note also that since
$\wkj(\cdot, t)$ lies in a compact subset of $C(\overline{\Omega})$ independently of $j$ and $t$
(since $\wkj (\cdot, t)$ lies in a bounded set in $C^{\gamma}(\Omega)$ for some $\gamma>0$, independent
of $j, t$), the fact that $\sup_{t \in \R} \| \wkj (\cdot, t) - \tilde{w} \|_{L^2(\Omega)} \rightarrow
0$ as $j \rightarrow \infty$ implies that $\sup_{t \in \R} \| \wkj (\cdot,t) -
\tilde{w} \|_{L^{\infty}(\Omega)} \rightarrow 0$ as $j \rightarrow \infty$.
This is because given a compact subset of $L^{\infty}(\Omega)$, the $L^{\infty}$-
and $L^2$-norms generate equivalent  metrics on this set due to the fact that a continuous
bijection on a compact set is a homeomorphism. 
Moreover, since $\Lambda \subset \subset \Omega$, it follows as in the proof of Lemma
\ref{convlem} (i) that
\begin{equation}
\label{cptconv}
~~~~\| \alpha \ukj (\cdot, t) - (\wkj)^+(\cdot, t) \|_{L^{\infty}(\tilde{\Lambda} \times [T, \infty))} \rightarrow 0
~~\mbox{and}~~\| \vkj (\cdot, t) + (\wkj)^-(\cdot, t) \|_{L^{\infty}(\tilde{\Lambda} \times [T, \infty))} \rightarrow 0
,
\end{equation}
as $j \rightarrow \infty$ for each fixed $T \in \R$ and $\tilde{\Lambda}$ a compact subset of $\Omega$.
So taking $\tilde{\Lambda}$ with $\Lambda \subset \subset \tilde{\Lambda}  \subset \subset  \Omega$, we have
the existence of $j_0$ such that  $\frac{x'}{\sqrt{\kj}}+\xj \in \tilde{\Lambda}$ for $j \geq j_0$ for all
$x' \in K \subset \RN$ compact (where $j_0$ is independent of $x'$ for a given $K$), and hence the uniform
convergence in (\ref{cptconv}) gives the existence of $\xb \in \Lambda$ such that
$$ \vkj (\frac{x'}{\sqrt{\kj}}+\xj, \frac{t'}{\kj}) \rightarrow - \tilde{w}^-(\xb), ~~
\alpha \ukj (\frac{x'}{\sqrt{\kj}}+\xj, \frac{t'}{\kj}) \rightarrow  \tilde{w}^+(\xb),
$$
for a subsequence 
as $j \rightarrow \infty$, uniformly in $(x', t') \in K \times [-T, T]$ for every $T>0$
and $K \subset \RN$ compact. We thus obtain an $L^{\infty}$-solution $(\tilde{h}, \tilde{l})$
on $\RN \times \R$ of
\begin{eqnarray}
h_t & = & \Delta h + \tilde{w}^-(\xb)h - \alpha^{-1} \tilde{w}^+(\xb) l, \label{limithl} \\
l_t & = & \Delta l -  \tilde{w}^+(\xb)l + \alpha \tilde{w}^- (\xb)h, \nonumber
\end{eqnarray}
such that $|\tilde{h}(0,0)| \geq \epz$ and $|\tilde{l}(0,0)| \geq \epz$.

Now note that since $\xb \in \Lambda$ and $\Lambda$ is compactly contained in $\Omega^0 \setminus
\{ x: \tilde{w} (x) =0 \}$, exactly one of $\tilde{w}^+(\xb)$ and  $\tilde{w}^-(\xb)$ is non-zero.
Suppose $\tilde{w}^+(\xb) \neq 0$. Then $\tilde{w}^+(\xb) >0$ and
\begin{equation}
\tilde{l}_t = \Delta \tilde{l} -   \tilde{w}^+(\xb) \tilde{l} ~~~\mbox{for all}~~(x,t) \in \RN \times \R.
\label{ltilde}
\end{equation}
If $\sup_{(x,t) \in \RN \times \R} \tilde{l} (x,t) = \sup_{(x,t) \in \RN \times \R} \{ - \tilde{l} \xt \}=0$,
then
$\tilde{l}(x,t)=0$ for all $(x,t)$, which contradicts $|\tilde{l}(0,0)| \geq \epz>0$. Otherwise,
either $\sup_{(x,t) \in \RN \times \R} \tilde{l} (x,t) >0$ or $\sup_{(x,t) \in \RN \times \R}\{ - \tilde{l}
(x,t)\} >0$; in the latter case, replace $\tilde{l}$ by $-\tilde{l}$ (which still satisfies (\ref{ltilde}))
and $|- \tilde{l}(0,0)|>0$. Now as in the proof of Lemma \ref{keylem}, define
$$z(t) = \sup_{x \in \RN} \tilde{l} (x,t). $$
Arguing as in the proof of Lemma \ref{keylem} then gives that
\begin{equation}
\dot{z}(t) \leq -  \tilde{w}^+(\xb) z(t)~~~a.e.~ t \in \R.
\label{z2eqn}
\end{equation}
Now since $\sup_{(x,t) \in \RN \times \R} \tilde{l} (x,t) >0$, there exists $t$ with $z(t)>0$, so
$z(s) \geq z(t)>0$ for all $s \leq t$, since (\ref{z2eqn}) implies that $z$ is non-increasing when
it is non-negative. Hence for any $t_0 < t \in \R$,
\begin{equation}
\label{z3eqn}
z(t) \leq z(t_0) \exp(- \tilde{w}^+(\xb)(t-t_0)),
\end{equation}
and so since $t_0 < t$ was arbitrary, we can let $t_0 \rightarrow -\infty$ in
(\ref{z3eqn}) to find that $z(t) \leq 0$, which contradicts the above. Similarly,
if $\tilde{w}^-(\xb) \neq 0$, the equation for $\tilde{h}$ yields a contradiction.
Hence the claim (\ref{hkjlkjsmall}) is true.

It remains to establish that $\hat{w}^{k_j}:= \alpha \hkj - \lkj$ is uniformly small on
$\clo \times \R$ if $j$ is large. The argument given for the corresponding result
in \cite[Thm 3]{DZ} applies almost un-changed and we omit the details.
Note that the requirement that $\tw$ be a non-degenerate solution of $\stat$ is needed here.
The idea is that via a contradiction argument, a non-trivial bounded solution of
the linearisation of the parabolic equation satisfied by the limit $\wh$ as $k_j \rightarrow \infty$
of $\hat{w}^{k_j}$ is obtained. Since this solution is non-trivial, there exists a time $\bar{t}$
such that $\wh (\cdot, \bar{t}) \not\equiv 0$, and by the non-degeneracy assumption, there must
be a {\em non-zero} real eigenvalue $\lambda$ of the linearisation of $\stat$ such that the
$L^2$-inner product of $\wh (\cdot, \bar{t})$
with a corresponding eigenfunction $\phi$ is non-zero. But then $z(t) := < \wh \fixt, \phi >$
can be shown to satisfy $\dot{z} = \lambda z$, and thus cannot be bounded, which is a contradiction.
Note that
having $\hat{w}$ satisfy zero Dirichlet rather than zero Neumann conditions causes
no difficulties, and that
the fact that
$$ \| \alpha \ukj (\cdot, t) - \tilde{w}^+\|_{L^{\infty}(\Lambda \times [-T, T])} \rightarrow 0
~~\mbox{and}~~  \| \vkj (\cdot, t) + \tilde{w}^-\|_{L^{\infty}(\Lambda \times [-T, T])} \rightarrow 0$$
as $j \rightarrow \infty$ for each $T>0$ and each $\Lambda \subset \subset \Omega$,
is enough to pass to the limit in the various weak forms of equations obtained by multiplying
by smooth functions of compact support in $\Omega \times \R$ (see \cite[Thm 3]{DZ}).

To conclude, suppose that $\Lambda$ is a compact subset of $\Omega \setminus \{x: \tilde{w}(x)=0 \}$.
Since $\hat{\wkj} = \alpha \hkj - \lkj$ converges uniformly to zero on $\Omega \times \R$ and
since, by (\ref{hkjlkjsmall}),  given $\epsilon >0$ there exists $j_0$ such that $j \geq j_0$ implies
$|\hkj (x,t)| < \epsilon$ or $|\lkj (x,t)|< \epsilon$ for each $x \in \Lambda, t \in \R$
(by the blow-up argument above), it follows that $\lkj$ and $\hkj$ each converge uniformly
to zero on $\Lambda \times \R$ as $j \rightarrow \infty$. Hence given $\hat{\epsilon} >0$, there
exists $\hat{j}$, independent of $\hat{t}$, such that for all $j \geq \hat{j}$,
$$ \int_{\Omega \times [\hat{t},\hat{t}+1]} (\hkj)^2 ~~ \leq ~~  \hat{\epsilon} + \int_{(\Omega \setminus \Lambda)
\times [\hat{t}, \hat{t}+1]} ( \hkj )^2 ~~~ \leq ~~~ \hat{\epsilon} +
( \| \hkj \|_{L^{\infty}(\Omega \times \R)})^2 | \Omega \setminus \Lambda |. $$
Now $\| \hkj \|_{L^{\infty}(\Omega \times \R)}$ is bounded independently of $j$ and $| \Omega \setminus \Lambda |$
can be made arbitrarily small by a suitable choice of $\Lambda$, since $\{ \: \tilde{w}(x)=0 \}$ and
$\partial \Omega$ each have zero $n$-dimensional measure. So there exists $\tilde{j}$ such that for all
$t \in \R$, $j \geq \tilde{j}$ implies that $\| \hkj \|_{L^2(\Omega \times [t, t+1])} \leq 1/8$.
A similar estimate for $\| \lkj \|_{L^2(\Omega \times [t, t+1])}$ can be established, giving a
contradiction with the normalisation (\ref{normalise}) for $\| (\hkj, \lkj) \|'$. The result follows. \qed

\begin{theorem}
\label{thmconvsyst}
Suppose that a solution $\tilde{w}$ of $\stat$ is non-degenerate. Then given $M>0$, there exist
$\epsilon, \kz >0$ such that if $ k \geq \kz$ and the solution $(\uk, \vk)$ of $\pk$ satisfies
$0 \leq \uk, \vk \leq M$ and
\begin{equation}
\label{closeconvsyst}
\| \wk (\cdot, t) - \tilde{w} \|_{L^2(\Omega)} \leq \epsilon
\end{equation}
for all $t$ sufficiently large, then there exists a non-negative stationary solution $(\tilde{\uk},
\tilde{\vk})$ of $\pk$ such that $\uk (\cdot, t) \rightarrow \tilde{\uk} (\cdot)$
and $\vk (\cdot, t) \rightarrow \tilde{\vk} (\cdot)$
in $\wtp$ and in $C^{1, \lambda'} (\overline{\Omega})$
for all $\lambda' \in (0, \lambda)$ as $t \rightarrow \infty$.
\end{theorem}

\pf Let $\epsilon, \kz$ be as in Lemma \ref{lemomegalimit}, and for (fixed) $k \geq \kz$, let
$(\uk, \vk)$ satisfy the hypotheses above (that such $(\uk, \vk)$ exist if there
are  non-degenerate solutions of $\stat$
follows from Theorem \ref{energythm}).
Let $\Gamma$ denote the $\omega$-limit set of $(
u_0^k , v_0^k)$ in $\wtp \times \wtp$. Recall Lemma \ref{lemregandcompact} (d) and note that the fact that
$\Gamma$ is invariant implies that it
consists of the union of trajectories of $\pk$ that
are defined for all $t \in \R$. Now it follows from (\ref{closeconvsyst}) that each
$(\gamma_u, \gamma_v) \in \Gamma$ satisfies
$\| \alpha \gamma_u -  \gamma_v - \tilde{w} \|_{L^2(\Omega)} \leq \epsilon.$
And by the characterisation of the omega-limit set, given $(\tilde{\gamma}_u, \tilde{\gamma}_v) \in \Gamma$,
there exists a solution $(\eta_u, \eta_v)$ of $\pk$, defined for all $t \in \R$, such that
\begin{description}
\item[(i)] $(\eta_u (\cdot, t), \eta_v (\cdot, t)) \in \Gamma$ for every $t \in \R$, and
\item[(ii)] $(\tilde{\gamma}_u, \tilde{\gamma}_v) = (\eta_u (\cdot, \hat{t}), \eta_v (\cdot, \hat{t}))$
for some $\hat{t} \in \R$.
\end{description}
So $\eta := \alpha \eta_u - \eta_v$ satisfies $\| \eta (\cdot, t) - \tilde{w} \|_{L^2(\Omega)} \leq \epsilon$
for all $t \in \R$, and thus it follows from Lemma \ref{lemomegalimit} that $(\eta_u, \eta_v)$ must be
independent of time; that is, $(\tilde{\gamma}_u, \tilde{\gamma}_v) = (\eta_u (\cdot, t), \eta_v (\cdot, t))$
for all $t \in \R$ and is a stationary solution of $\pk$. Hence $\Gamma$ consists entirely of stationary
solutions of $\pk$ (and, since $\Gamma$ is non-empty by Lemma \ref{lemregandcompact}, such solutions
must exist).

Since dist$(\ukvk \fixt, \Gamma) \rightarrow 0$ as $t \rightarrow \infty$, by Lemma \ref{lemregandcompact} (d),
it remains to show that for $k$ sufficiently large,
the elements of $\Gamma$ are isolated in $\wtp \times \wtp$.
Consider $F: \wtp \times \wtp \rightarrow L^p (\Omega) \times  L^p (\Omega)$ defined by
$$ F(u,v) = \left\{ \begin{array}{l} \Delta u + f(u) - kuv \\ \Delta v + g(v) - \alpha k u v. \end{array} \right.$$
Then since $\wtp \hookrightarrow C^{0, \lambda}(\Omega)$ for some $\lambda >0$,
$F \in C^1(\wtp \times \wtp,
L^p (\Omega) \times  L^p (\Omega) )$. Moreover, an argument the same
as part of the proof of
\cite[Thm 1.2]{DG} gives that for $k$ sufficiently large, the Fr\'echet derivative
$F'(u^*, v^*)$ at a solution $(u^*, v^*)$ of $F(u,v)=0$ is injective on $\wtp \times \wtp$.
That $F'(u^*, v^*)$ is also surjective and has bounded inverse follows from the Fredholm
Alternative. The isolatedness of elements of $\Gamma$ is then
a consequence of the Inverse Function Theorem (see \cite[Thm 1.2]{AP}, for example) and
the result follows. \qed

\noindent We conclude with a result on simple dynamics for $\pk$ and some remarks.

\begin{theorem}
\label{thmsimpledynamics}
Suppose that {\em all} the solutions of $\stat$ are  non-degenerate. Then
there
exists $\kz >0$ such that if  $k \geq \kz$,
there exists a non-negative stationary solution
$(\tilde{\uk},
\tilde{\vk})$ of $\pk$ such that $\uk (\cdot, t) \rightarrow \tilde{\uk} (\cdot)$
and $\vk (\cdot, t) \rightarrow \tilde{\vk} (\cdot)$
in $\wtp$ and in $C^{1, \lambda'} (\overline{\Omega})$
for all $\lambda' \in (0, \lambda)$ as $t \rightarrow \infty$.
Note that $\kz$ is dependent on the boundary data $\mkone |_{\partial \Omega}, \mktwo
|_{\partial \Omega}$ but is independent
of the choice of initial data for $\pk$.
\end{theorem}

\pf Note first that there exists $M>0$  such that for any $k$ and any initial data $(u_0^k, v_0^k)$ for $\pk$
satisfying $(u_0^k, v_0^k) = (\mkone, \mktwo)$ on $\partial \Omega$
 where $\mkone, \mktwo$ satisfy $\bone - \bthree$, there exists $T$
(dependent on $k$ and $(u_0^k, v_0^k)$) such that $0 \leq \uk \fixt, \vk \fixt \leq M$ for all
$t \geq T$. This follows from the fact that $z(t) := \sup_{x \in \overline{\Omega} }
\uk (x,t)$ (similarly, $\sup_{x \in \overline{\Omega} } \vk (x,t)$) decreases at a certain rate
for $t$ for which $z(t) > \sup_{k \in \N, ~x \in \overline{\Omega}} \{ \mkone (x), \mktwo (x), 2 \}
=: M$; indeed one can check that $z_t
\leq f(z)$ so that $z$ lies below the solution of the ordinary
differential equation $U_t = f(U)$ together with the same initial
condition $ \sup_{x \in \bar \Omega} u_0(x)$. Thus $z (\bar{t}) <M$ for some $\bar{t}$; then, since $z(t)$ is decreasing whenever
$z(t) > \sup_{k \in \N, ~x \in \overline{\Omega}} \{ \mkone(x), \mktwo(x), 1 \}$, $z(t)$ cannot
increase above $M$ for any $t > \bar{t}$.

The result now follows  immediately from Theorem \ref{energythm} and Theorem \ref{thmconvsyst}
applied with this value of $M$. \qed

\section{On the local existence and uniqueness of stationary solutions of $\pk$ close to a non-degenerate solution of
$\stat$ for large $k$}
We first prove the following result on the existence (and total degree) of positive stationary
solutions of $\pk$ near $\splitstat$ in $\lplp$. Note that we assume here that the boundary
conditions $\mkone, \mktwo$ are in fact {\em independent} of $k$ and write $\mkone = m_1$,
$\mktwo = m_2$.

\begin{theorem}
\label{thmaone}
Suppose $\wz$ is an isolated (in $L^p(\Omega)$) solution of $\stat$ which changes sign and
has non-zero index. Suppose further that the boundary conditions $\mkone, \mktwo$ in
$\pk$ are independent of $k$. Then there exist $\kz$ and $\delta_1 >0$ such that
for $k \geq \kz$, $\pk$ has a positive stationary solution $(u,v)$ in the
$\delta_1$-neighbourhood in $\lplp$ of $\splitstat$. Here $p$ is as in condition
\bone and by the index of $\wz$ we mean the fixed point index
$$ \mbox{index}_K (B_2, \wz),$$
where $K = \{ w \in \coneob : w= \bc ~~\mbox{on}~~ \partial \Omega \}$ and
$B_2 w$ is the unique solution $y$ of
\begin{eqnarray}
\label{linearstat}
- \Delta y &=& \alpha f (\alpha^{-1} w^+) - g(- w^-) ~~~~~~~~~\mbox{in}~~\Omega, \\
y &=& \bc ~~~~~~~~~~~~~~~~ ~~~~~~~\mbox{on}~~ \partial \Omega. \nonumber
\end{eqnarray}
(Note that we use the notation $B_2$ here for ease of reference with \cite{DD}.)
\end{theorem}

\pf This is an analogue of \cite[Thm 3.3]{DD}, which establishes a similar result for a
system with homogeneous Dirichlet boundary conditions. Some key parts of the
proof differ from that of \cite[Thm 3.3]{DD} and so we include a proof here, giving
most detail where the differences lie. Consider the homotopy

\begin{equation}
\label{homo}
\begin{array}{ll}
- \Delta u = t f(u) + (1-t) f(\ucombp) - kuv & \mbox{in}~~ \Omega, \\
- \Delta v = t g(v) + (1-t) g(\vcombp) - \alpha kuv & \mbox{in}~~ \Omega, \\
u=m_1 & \mbox{on}~~\partial \Omega,\\
v=m_2 & \mbox{on}~~\partial \Omega,
\end{array}
\end{equation}
where $t \in [0,1]$.

We first note that positive solutions $(u,v)$ are bounded in $L^{\infty}(\Omega)$
independently of $t \in [0,1]$ and $k$. Indeed, it follows from {\bf (a)} that there exists
$c>0$, independent of $(u,v)$, $k$ and $t$, such that $- \Delta u \leq c$, $- \Delta v
\leq c$. Now let $y_1, y_2 \in \wtp$ be such that $- \Delta y_i = c$ in $\Omega$ and
$y_i = m_i$ on $\partial \Omega$, ($i=1,2$). Then the maximum principle gives that
$u \leq y_1$, $v \leq y_2$. Since $u,v \geq 0$, it follows that there exists a constant
$M_0 >0$ such that for any non-negative solution $(u,v)$ of (\ref{homo}),
\begin{equation}
\label{bounduv}
0 \leq u \leq M_0, ~~~~  0 \leq v \leq M_0.
\end{equation}
Now, as in \cite[Thm 3.3]{DD}, let $f_1(u,v,t)$ and $f_2(u,v,t)$ denote the right-hand-sides
of the equations for $u$ and $v$ respectively in (\ref{homo}), and define
$u_M = \min \{u, M \}$, $v_M = \min \{v, M \}$. Next define
$$ \tilde{f}_i (u,v,t) = f_i (u_{M_0+1}, v_{M_0+1}, t), ~~~i=1,2. $$
By the choice of $M_0$ in (\ref{bounduv}), the modified problem
\begin{equation}
\label{modprob}
\begin{array}{ll}
-\Delta u = \fit & \mbox{in}~~ \Omega, \\
-\Delta v = \ftt & \mbox{in}~~ \Omega, \\
u= m_1 & \mbox{on}~~ \partial \Omega, \\
v= m_2 & \mbox{on}~~ \partial \Omega,
\end{array}
\end{equation}
has the same non-negative solution set as (\ref{homo}).
Indeed every
nonnegative solution pair $(u,v)$ of (\ref{modprob}) is such that $u, v \leq
M_0$ so that $(u,v)$ satisfies (\ref{homo}).

Now choose $\delta >0$ small enough that  $\wz$ is the only solution of $\stat$
in the $\delta-$neighbourhood $\nw$ of $\wz$ in $L^p(\Omega)$. Then choose
$\delta_1 >0$ so that 
\begin{equation}
(u,v) \in \partial \npair ~\mbox{implies that}~~ u \not\equiv 0,~
v \not\equiv 0~~\mbox{and}~~\alpha u - v \in \nw.
\label{dhone}
\end{equation}
Here  $\partial \npair$ denotes the
boundary of the $\delta_1$-neighbourhood $\npair$ of $\splitstat$ in $\lplp$.\\ 

\noindent Next we prove the following result.  

\begin{lemma}
For the above choice of  $\done$, there exists $\kz$ such that (\ref{modprob})
has no non-negative solution $(u,v)$ with $(u,v) \in \partial \npair$ for any
$t \in [0,1]$ and $k \geq \kz$. 
\label{danlem}
\end{lemma}
\pf Because of our boundary conditions, more care is
needed here than in the corresponding
proof in \cite[Lemma 3.1]{DD}.
Suppose, for contradiction, that there are $k_n \rightarrow \infty$, $t_n \in [0,1]$
such that (\ref{modprob}) has a non-negative solution $(u_n, v_n) \in \partial \npair$.
Then by (\ref{bounduv}), $(u_n, v_n)$ is a solution of (\ref{homo}) for $k=k_n$,
$t=t_n$. Setting $w_n = \alpha u_n - v_n$ gives
\begin{equation}
\label{dhtwo}
~~~~ - \Delta w_n = \alpha (t_n f(u_n) + (1-t_n)f((u_n - \alpha^{-1} v_n)^+))
- t_n g(v_n) - (1-t_n)g((v_n - \alpha u_n)^+) =: b_n, ~\mbox{say},
\end{equation}
and $b_n$ is bounded in $L^{\infty}(\Omega)$ independently of $n$, again using
(\ref{bounduv}). Now let $\psi \in \wtp$ be such that $\Delta \psi =0$ in $\Omega$
and $\psi = \bc$ on $\partial \Omega$, and set $W_n = w_n - \psi$. Then
$- \Delta W_n = b_n$, and there is a constant $K>0$ such that
$$ \int_{\Omega}  | \nabla W_n |^2 = - \int_{\Omega} W_n \Delta W_n = \int_{\Omega} W_n b_n \leq K. $$
Hence $w_n$ is bounded in $W^{1,2}(\Omega)$ and thus, taking a subsequence if necessary, there
exists $\overline{w} \in W^{1,2}(\Omega)$ such that
\begin{equation}
\label{wnconv}
w_n \rightharpoonup \overline{w} ~~\mbox{in}~~ W^{1,2}(\Omega)~~\mbox{and}~~L^2(\partial \Omega),~~
\mbox{and}~~ w_n \rightarrow \overline{w} ~~\in~~ L^2 (\Omega).
\end{equation}
We now adapt an idea from \cite{CDHMN}. Let $\phi \in W^{1,2}_0(\Omega)$ satisfy
$- \Delta \phi = \lambda \phi$ in $\Omega$, $\phi =0$ on $\partial \Omega$ with
$\lambda >0$ and $\phi >0$ in $\Omega$. Then multiplication of (\ref{homo}) by $u_n \phi$, integration by parts and
the fact that $u_n, v_n \geq 0$ give
\begin{equation}
\label{unphi}
~~~~\int_{\Omega} | \nabla u_n |^2 \phi ~dx \leq \int_{\Omega} \{ t_n f(u_n)) + (1-t_n)f((u_n - \alpha^{-1}
v_n)^+) \} u_n \phi + \frac{1}{2} u_n^2 \Delta \phi ~dx - \int_{\partial \Omega} \frac{1}{2} u_n^2
\frac{\partial \phi}{\partial \nu} dS,
\end{equation}
and hence there exists $K_1>0$, independent of $n$, such that
\begin{equation}
\label{unbound}
\int_{\Omega} | \nabla u_n |^2 \phi~ dx \leq K_1.
\end{equation}
Similarly, there exists $K_2, K_3 >0$ such that
\begin{equation}
\label{vnbound}
\int_{\Omega} | \nabla v_n |^2 \phi~ dx \leq K_2,
\end{equation}
and
\begin{equation}
\label{unvn}
k_n \int_{\Omega} u_n v_n \phi ~ dx \leq K_3.
\end{equation}
(See \cite{CDHMN} for details of similar arguments for the parabolic system $\pk$.)
It follows from (\ref{unbound}), (\ref{vnbound}) and (\ref{bounduv}) that $\{ u_n \}_{n=1}^{\infty}$,
$\{ v_n \}_{n=1}^{\infty}$ are each bounded in $W^{1,2} (\Omega ')$ and hence relatively compact in
$L^2(\Omega')$ for each $\Omega' \subset \subset \Omega$. Using (\ref{bounduv}) again, 
it follows that there are subsequences of $u_n, v_n$ (not relabelled) and
$\overline{u}, \overline{v} \in L^{\infty}(\Omega)$ such that
\begin{equation}
\label{compactp}
u_n \rightarrow \overline{u}, ~~v_n \rightarrow \overline{v} ~~\mbox{in}~~L^p(\Omega)~~\mbox{and}~~a.e.~\mbox{in}
~~\Omega,
\end{equation}
and by (\ref{unvn}),
\begin{equation}
\label{seg}
\overline{u} ~\overline{v} = 0 ~~a.e. ~~\mbox{in}~~ \Omega,
\end{equation}
which, together with (\ref{wnconv}), gives that $\overline{u} = \alpha^{-1} \overline{w}^+$ and
$\overline{v} = - \overline{w}^-$. Note that (\ref{compactp}) also gives that $(\overline{u}, \overline{v})
\in \partial \npair$.

Now take a subsequence if necessary to ensure $t_n \rightarrow \overline{t} \in [0,1]$.
It follows from (\ref{dhtwo}), (\ref{wnconv}) and (\ref{compactp}) 
that for $\phi \in C^{\infty}_0(\Omega)$,
$$\int_{\Omega} \nabla \overline{w} \nabla \phi ~ dx = \int_{\Omega} \{\alpha (\overline{t} f(\alpha^{-1}\overline{w}^+)
+ (1- \overline{t}) f(\alpha^{-1} \overline{w}^+)) - \overline{t} g(- \overline{w}^-) -
(1- \overline{t}) g(- \overline{w}^-)\} \phi ~dx. $$
Thus $\overline{w}$ is a solution of $\stat$. By (\ref{dhone}), (\ref{wnconv}) and
(\ref{compactp}), $\overline{w} = \alpha
\overline{u} - \overline{v} \in \nw$, and thus $\overline{w} = \wz$, so $\overline{u} = \alpha^{-1}
\wz^+$ and $\overline{v} = - \wz^-$. This contradicts that $(\overline{u}, \overline{v}) \in
\partial \npair$, and completes the proof of Lemma \ref{danlem}. \qed

\noindent  Next we return to the proof of Theorem \ref{thmaone}. 
Now given $k \geq \kz$, choose $M_k >0$ sufficiently large that
\begin{equation}
\label{fitmk}
\fit + M_k u \geq 0, ~~\ftt + M_k v \geq 0
\end{equation}
for any $u,v \geq 0$ and $t \in [0,1]$, and also
\begin{equation}
\label{fitderiv}
\frac{\partial}{\partial u} \fit + M_k >0, ~~~\frac{\partial}{\partial v} \ftt + M_k >0
\end{equation}
for $0 \leq u,v \leq M_0$ and $t \in [0,1]$ (note that $\fit, \ftt$ are differentiable 
at such $u,v$). Define
\begin{equation}
\label{compopdef}
A_t = A_{t,k}: \lplp \rightarrow \lplp
\end{equation}
by $A_t(u,v) = (y,z)$, where
\begin{equation}
\label{danlabel}
 \begin{array}{ll}
(- \Delta + M_k) y = \fit + M_k u & \mbox{in} ~\Omega, \\
(- \Delta + M_k) z = \ftt + M_k v & \mbox{in} ~\Omega, \\
y=m_1 & \mbox{on}~~ \partial \Omega,\\
z=m_2 & \mbox{on}~~ \partial \Omega.\\
\end{array} \end{equation}
Then $A_t$ is completely continuous (that is, $A_t$ is continuous and compact) and maps the natural positive cone $P$ in $\lplp$ into
itself. Moreover, by  Lemma \ref{danlem} and the homotopy invariance of the degree
(see, for example, \cite[p201]{AD}), for
$k \geq \kz$,
\begin{equation}
\label{degeq}
\mbox{deg}_P(I-A_0, P \cap \npair, 0) = \mbox{deg}_P (I- A_1, P \cap \npair, 0).
\end{equation}
Note that $(u,v) = A_0 (u,v)$ if and only if $(u,v)$ solves (\ref{modprob}) with
$t=0$, and by (\ref{bounduv}), such $(u,v)$ satisfies
\begin{equation}
\label{endpoint}
\begin{array}{ll}
- \Delta u = f(\ucombp) - kuv & \mbox{in}~~ \Omega, \\
- \Delta v = g(\vcombp) - \alpha kuv & \mbox{in} ~~\Omega,\\
u=m_1 & \mbox{on}~~\partial \Omega, \\
v=m_2 & \mbox{on}~~\partial \Omega.
\end{array}
\end{equation}
As in \cite[Thm 3.3]{DD}, we will show that for $k$ sufficiently large, (\ref{endpoint})
has a unique non-negative solution in $\npair$. Note first that if $(u,v) \in \npair$ is a non-negative
solution of (\ref{endpoint}), then $\tilde{w}_0 := \alpha u - v$ is a solution of $\stat$. And hence
$\tilde{w}_0 = \wz$, by the choice of $\delta_1$. It follows that $\alpha u - v= \wz$, and, since $u,v \geq 0$,
that
\begin{equation}
\label{ulowerbound}
u \geq \alpha^{-1} \wz^+.
\end{equation}

\noindent Now observe that for any $k>0$, the equation
\begin{equation}
\label{oneeq}
- \Delta u = f( \alpha^{-1} \wz^+) - ku(\alpha u - \wz)~~\mbox{in}~~\Omega,~~ u=m_1~~\mbox{on} ~\partial \Omega,
\end{equation}
has a unique solution $u^k$ satisfying $\uk \geq \alpha^{-1} \wz^+$. Indeed, $\alpha^{-1} \wz^+$
is a lower solution of (\ref{oneeq}). This is because Kato's inequality (see, for example, \cite{Kato, DS}) gives that in the
sense of distributions, 
\begin{eqnarray}
- \Delta (|\wz|) & \leq & - \mbox{sign}( \wz)~ \Delta \wz \nonumber \\
& = & \mbox{sign}( \wz)~ (\alpha f(\alpha^{-1} \wz^+) - g(-\wz^{-})), ~~\mbox{(since $\wz$ satisfies (\ref{stationary}))}, \nonumber 
\\ & = & \alpha f(\alpha^{-1} \wz^+) + g(-\wz^{-}). \label{kato}
\end{eqnarray}
Then because $|\wz|+\wz = 2\wz^+$, adding (\ref{kato}) and (\ref{stationary}) gives that in the sense of
distributions, 
$$ - \Delta (\alpha^{-1} \wz^+) \leq f(\alpha^{-1} \wz^+) = f(\alpha^{-1} \wz^+)
-k \alpha^{-1}\wz^+(\alpha (\alpha^{-1} \wz^+) - \wz), $$
and thus since $\alpha^{-1} \wz = m_1 - \alpha^{-1}m_2 \leq m_1$ on $\partial \Omega$ and hence $\alpha^{-1} \wz^+ \leq m_1$ on $\partial \Omega$, it follows that $\alpha^{-1} \wz^+$ is a lower solution of (\ref{oneeq}). 
And any large positive constant is an upper solution of (\ref{oneeq}), so there
exists a solution $\uk \geq \alpha^{-1} \wz^+$. Uniqueness follows from the fact that the right-hand-side
of (\ref{oneeq}) is non-increasing in $u$ when $u \geq \alpha^{-1} \wz^+$. (Note that this argument
differs from that used in \cite[Thm 3.3]{DD} since here $u \equiv 0$ is no longer necessarily a lower
solution if $\alpha m_1 - m_2$ is large on $\partial \Omega$. )

Next, we prove that $\uk \rightarrow \alpha^{-1} \wz^+$ in $L^p (\Omega)$ as $k \rightarrow \infty$.
This argument is the same as that in \cite[Thm 3.3]{DD}; we give it here for completeness.
Note that $\uk \geq \alpha^{-1} \wz^+$ and $\uk (\alpha \uk - \wz) \geq 0$. This second inequality
gives that $u^{k_1}$ is an upper solution of (\ref{oneeq}) if $k_1 \leq k$ and hence $u^{k_1} \geq \uk$
if $k_1 \leq k$. Now let $\lim_{k \rightarrow \infty} \uk (x) = \overline{u}(x), ~~x \in \Omega$. Then
for any $\phi \in C^{\infty}_0(\Omega)$,
$$ \int_{\Omega} \uk (\alpha \uk - \wz) \phi = k^{-1} \left( \int_{\Omega} \uk \Delta \phi + \int_{\Omega}
f(\alpha^{-1} \wz^+) \phi \right) \rightarrow 0$$
as $ k \rightarrow \infty$, and hence $$\int_{\Omega} \overline{u} (\alpha \overline{u} - \wz ) \phi =0~~~
\mbox{for all}~~\phi \in C^{\infty}_0 (\Omega). $$
So $ \overline{u} (\alpha \overline{u} - \wz )=0$. Hence $\overline{u} = \alpha^{-1} \wz^+$, since
we know $\overline{u} \geq \alpha^{-1} \wz^+$. Now let $\vk = \alpha \uk - \wz$. Then $(u,v) = (\uk, \vk)$
solves (\ref{endpoint}) and
$ (\uk, \vk) \rightarrow (\alpha^{-1} \wz^+, - \wz^-)$ in $\lplp$ as $ k \rightarrow \infty$.
Thus $(\uk, \vk)$ is a non-negative solution of (\ref{endpoint}) in $\npair$ when $k$ is sufficiently large.

Conversely, if $(u,v) \in \npair$ is a non-negative solution of (\ref{endpoint}), then $\overline{w}_0
= \alpha u - v$
solves $\stat$ and $\overline{w}_0 \in \nw$ so $\overline{w}_0 = \wz$. Thus
$$ - \Delta u = f(\alpha^{-1} \wz^+) - ku(\alpha u - \wz)~~\mbox{in}~~\Omega, ~~~u=m_1~~\mbox{on}~~\partial
\Omega,$$ and hence $u=\uk$ and $v=\alpha u - \wz = \vk$.

Thus for $k$ sufficiently large, $(\uk, \vk)$ is the unique non-negative solution of (\ref{endpoint})
in $\npair$. So there exists $k_1 \geq \kz$ such that for $k \geq k_1$,
\begin{equation}
\label{deg1}
\mbox{deg}_P(I- A_0, P \cap \npair, 0) = \mbox{index}_P (A_0, (\uk, \vk)).
\end{equation}
Let
\begin{equation}
\label{setC}
C= \{ (u,v) \in \coneob \times \coneob : u,v \geq 0 ~~\mbox{in}~~\Omega~~\mbox{and}~~u=m_1,~v=m_2~\mbox{on}~
\partial \Omega \}.
\end{equation}
Since $A_0$ maps $P$ into $C$,
it follows from two applications of the commutativity of the fixed point index (\cite[p 214]{Deimling},
\cite{Nussbaum}) that
\begin{equation}
\label{indexeq}
\mbox{index}_P(A_0, (\uk, \vk)) = \mbox{index}_C (A_0, (\uk, \vk)) = \mbox{index}_{\tilde{C}} (\tilde{A_0},
(\uk - h_1, \vk - h_2)),
\end{equation}
where $\tilde{C} = \{ (\tilde{u}, \tilde{v}) \in C^1_0 (\overline{\Omega}) \times C^1_0 (\overline{\Omega}):
\tilde{u} + h_1 \geq 0, ~ \tilde{v} + h_2 \geq 0 \}$, $h_i, ~i=1,2$ satisfy $\Delta h_i =0$ in $\Omega$
and $h_i = m_i$ on $\partial \Omega$, and for $(\tilde{u}, \tilde{v}) \in C^1_0 (\overline{\Omega}) \times
C^1_0 (\overline{\Omega})$,
\begin{equation}
 \tilde{A}_0 (\tilde{u}, \tilde{v}) = (- \Delta + M_k)^{-1}( \tilde{f}_1 (\tilde{u}+h_1, \tilde{v}+h_2, 0)
+ M_k \tilde{u}, \tilde{f}_2 (\tilde{u}+h_1, \tilde{v}+h_2, 0) + M_k \tilde{v}) 
\label{newop}
\end{equation}
where  the inverse $(-\Delta + M_k)^{-1}$ is taken under zero Dirichlet boundary conditions.\\

Now the strong maximum principle gives that $\uk>0$, $\vk >0$ in $\Omega$, and that the outward normal
derivative $\displaystyle{\frac{\partial \uk (\vk)}{\partial \nu} <0}$ at a point on $\partial \Omega$ where $\uk (\vk) =0$,
and hence $(\uk-h_1, \vk-h_2) \in \mbox{int} ~\tilde{C}$. So
\begin{equation}
\label{index2}
\mbox{index}_{\tilde{C}}(\tilde{A_0},
(\uk - h_1, \vk - h_2)) = \mbox{index}_{\conez} (\tilde{A_0},
(\uk - h_1, \vk - h_2)).
\end{equation}
(Note that our inhomogeneous boundary conditions necessitate a slightly different argument here from
that in \cite[Thm 3.3]{DD}.)

We now use the homeomorphism $h(u,v) = (u, \alpha u - v)$ in $\conez$, as in \cite[Thm 3.3]{DD}.
Note that $h^{-1} =h$, and that for $(\tilde{u}, \tilde{w}) \in \conez$,
\begin{small}
\begin{eqnarray*}
\lefteqn{(h^{-1} \tilde{A_0} h)(\tilde{u}, \tilde{w})  =  h^{-1}  \tilde{A_0} (\tilde{u}, \alpha \tilde{u} - \tilde{w})} \\
& = & h^{-1}  (- \Delta + M_k)^{-1}  \\ \lefteqn{ \left( \begin{array}{l}
f(\alpha^{-1}(\tilde{w} + \alpha h_1 - h_2)^+) - k(\tilde{u}+h_1)(\alpha(\tilde{u}+h_1)-(\tilde{w}+\alpha h_1 - h_2))
+M_k \tilde{u} \\
g((-(\tilde{w} + \alpha h_1-h_2))^-) - \alpha k (\tilde{u}+h_1)(\alpha(\tilde{u}+h_1)-(\tilde{w}+\alpha h_1 - h_2))
+ M_k (\alpha \tilde{u} - \tilde{w}) \end{array} \right)  } \\
& = & (- \Delta + M_k)^{-1} \\ & & \left( \begin{array}{l} f(\alpha^{-1}(\tilde{w} + \alpha h_1 - h_2)^+) -
k(\tilde{u}+h_1)(\alpha(\tilde{u}+h_1)-(\tilde{w}+\alpha h_1 - h_2))
+M_k \tilde{u} \\
\alpha f(\alpha^{-1}(\tilde{w} + \alpha h_1 - h_2)^+) - g(-(\tilde{w} + \alpha h_1-h_2)^-) + M_k\tilde{w}  \end{array} \right).
\end{eqnarray*}
\end{small}
So the commutativity of the fixed point index and the product formula give that
\begin{equation}
\label{productform} \mbox{index}_{\conez}(\tilde{A_0},(\uk - h_1, \vk - h_2)) = \mbox{index}_{C^1_0(\overline{\Omega})}
(\tilde{B}_2, \wz - \alpha h_1+h_2). \mbox{index}_{C^1_0(\overline{\Omega})} (\tilde{B}, \uk - h_1) 
\end{equation}
where
\begin{equation}
\label{btwotilde}
 \tilde{B}_2 \tilde{w} = (- \Delta + M_k)^{-1} ( \alpha f(\alpha^{-1}(\tilde{w} + \alpha h_1 - h_2)^+) - g(-(\tilde{w} + \alpha h_1-h_2)^-) + M_k \tilde{w} ) 
\end{equation}
and
$$ \tilde{B} \tilde{u} = (- \Delta + M_k)^{-1}  (f(\alpha^{-1}w_0^+) -
k(\tilde{u}+h_1)(\alpha(\tilde{u}+h_1)-w_0)
+M_k \tilde{u}) $$
for $(\tilde{u}, \tilde{w}) \in \conez$. Hence, by (\ref{indexeq}) and (\ref{productform}),
\begin{equation}
\label{index3}
\mbox{index}_P(A_0, (\uk, \vk)) = \mbox{index}_K (B_2, \wz) .\mbox{index}_{C^1_0(\overline{\Omega})}
(\tilde{B}, \uk - h_1)
\end{equation}
where $B_2, K$ are as defined in the statement of the theorem.

It remains to show that $\mbox{index}_{C^1_0(\overline{\Omega})}(\tilde{B}, \uk-h_1)=1$.
First note that $\alpha^{-1} \wz^+$ and $u_0$ (the solution of (\ref{oneeq}) with $k=0$)
are lower and upper solutions for (\ref{oneeq}) respectively (neither of which are solutions).
Just as in the system (\ref{danlabel}) we add to
                 both sides of (\ref{oneeq}) a term of the form $M_k u$ with
                 $M_k$ large enough so that $\tilde B$ maps the set
$$ C^* =  \{ \tilde{u} \in  C^1_0 (\overline{\Omega}) : \alpha^{-1} \wz^+ - h_1 \leq \tilde{u} \leq u_0 - h_1 \} $$
into itself (in fact, into $\mbox{int} C^*$, by the strong maximum principle).
Also, by the uniqueness for (\ref{oneeq})
discussed above, $\uk - h_1$ is the only solution of $\tilde{u}=\tilde{B}\tilde{u}$ in $C^*$. Also, $\tilde{B}(C^*)$
is a bounded set in $C^1_0(\overline{\Omega})$; choose a large ball $B_R$ in $C^1_0(\overline{\Omega})$
such that $\tilde{B}(C^*) \subset B_R$. Let $S=C^* \cap B_R$. Then $S$ is a bounded convex set in $C^1_0(\overline{\Omega})$
and $\tilde{B}$ maps $S$ into itself. Moreover, $\tilde{B}$ has a unique fixed point in int$S$. Thus
\begin{equation}
\label{index4}
\mbox{index}_{C^1_0(\overline{\Omega})} (\tilde{B}, \uk - h_1) = \mbox{deg}(I-\tilde{B}, \mbox{int}S, 0) =1.
\end{equation}
It follows from (\ref{degeq}), (\ref{deg1}), (\ref{index3}), (\ref{index4}) and the hypotheses of Theorem \ref{thmaone} 
that
\begin{equation}
\label{finaldeg}
\mbox{deg}_P(I-A_1, P \cap \npair, 0) = \mbox{index}_K(B_2, \wz) \neq 0,
\end{equation}
from which the result follows by the existence property of degree. \qed

\begin{theorem}
\label{thmuniqueness}
Suppose $\wz$ is a non-degenerate solution of $\stat$. Then there exists $\kz >0$ such that
for any $k \geq \kz$, $\pk$ has a unique positive stationary solution $(u,v)$ near
$\splitstat$ in $\lplp$.
\end{theorem}

\pf We mimic the proof of \cite[Thm 1.2]{DG} which treats the corresponding problem with
homogeneous rather than inhomogeneous Dirichlet boundary conditions. That in turn draws on
\cite[Thm 3.3]{DD}, of which our analogue is Theorem \ref{thmaone} above.

First, it follows exactly as in the proof of \cite[Thm 1.2]{DG} that if $\delta_1>0$ is sufficiently small, then for large enough $k$, the linearisation of the
stationary system
$$ \pks ~~~~\begin{array}{ll}
- \Delta u = f(u) - kuv & \mbox{in}~~ \Omega,\\
- \Delta v = g(v) - \alpha kuv & \mbox{in}~~ \Omega,\\
u=m_1 & \mbox{on}~~\partial \Omega,\\
v=m_2 & \mbox{on}~~\partial \Omega,
\end{array}
$$
is invertible at  positive solutions $(\uk, \vk) \in \npair $   (note that $p>N$, from \bone).

The next step is to show that for each such $k$, the value of $\mbox{index}_P(A_{1,k}, (\uk, \vk))$ at any positive solution $(\uk, \vk)$ of $\pks$ in
$\npair$  is the same, where
$P$ is (as above) the natural positive cone in $\lplp$ and $A_{1,k}$ is defined in
(\ref{compopdef}) (using (\ref{fitmk}) and (\ref{fitderiv})) in the proof of Theorem \ref{thmaone}.
This, together with Theorem \ref{thmaone}, will give the existence of a unique positive solution of
$\pks$ in $\npair$, since Theorem \ref{thmaone} gives the existence of such solutions, and the proof of Theorem
\ref{thmaone} gives that if $\delta_1>0$ is sufficiently small, then for large enough $k$, 
\begin{equation}
\label{degindexfinal}
\mbox{deg}_P(I-A_{1,k}, P \cap \npair, 0) = \mbox{index}_K(B_2, \wz),
\end{equation}
(see (\ref{finaldeg})). Note that $\mbox{index}_K(B_2, \wz)$ is $+1$ or $-1$
since $\wz$ is a non-degenerate fixed point of $B_2$ ($B_2$ is defined in (\ref{linearstat}) and in (\ref{btwotilde})).
Now for $(u^k, v^k)$ as above, it follows from the commutativity of the fixed point index that just as in (\ref{indexeq}),
$$ \mbox{index}_P(A_{1,k},(\uk, \vk)) = \mbox{index}_C (A_{1,k},(\uk, \vk))$$
for positive solutions of $\pks$
where $C$ is as defined in (\ref{setC}). And note (as in \cite{DG}), that in $C$, small
neighbourhoods of solutions are uniformly close, and hence the truncations in the definition of
$A_{1,k}$ do not affect $A_{1,k}$ near fixed points. Thus we can delete the truncations and work
with the map $\hat{A}_{1,k}$ defined to be $A_{1,k} |_C$ without the truncations. As in \cite{DG},
the reason for doing this is that $\hat{A}_{1,k}$ is differentiable. Now (see also (\ref{indexeq}) and (\ref{newop}))
\begin{eqnarray*}
\mbox{index}_C(\hat{A}_{1,k}, (\uk, \vk)) & = & \mbox{index}_{\tilde{C}}
(\widetilde{\hat{A}_{1,k}}, (\uk-h_1,\vk-h_2)) \\
& & \mbox{where}~~ \widetilde{}~~ \mbox{has an analogous effect to that in (\ref{newop}) }\\
& = & \mbox{index}_{\conez}(\widetilde{\hat{A}_{1,k}},  (\uk-h_1, \vk-h_2)) \\
& = & \mbox{index}_{\conez}(\widetilde{\hat{A}_{1,k}}'(\uk-h_1, \vk-h_2), 0)\\
& = & \mbox{index}_{\conez}(\hat{A}_{1,k}'(\uk, \vk), 0).
\end{eqnarray*}
Further,  it follows from the proof of Theorem \ref{thmaone} that the system
\begin{equation}
\label{hybridsyst}
\begin{array}{ll}
- \Delta u = f(\ucombp) - kuv & \mbox{in}~~\Omega,\\
- \Delta v = g(\vcombp) - \alpha kuv & \mbox{in}~~\Omega,\\
u=m_1 & \mbox{on}~~\partial \Omega,\\
v=m_2 & \mbox{on}~~\partial \Omega,
\end{array}
\end{equation}
has a unique solution $(\overline{u}_k, \overline{v}_k )$ such that $(\overline{u}_k, \overline{v}_k )
\rightarrow \splitstat$ in $\lplp$. Note that solutions of
(\ref{hybridsyst}) are fixed points of the operator $A_{0,k}$ from the proof of Theorem \ref{thmaone}.
And it follows as in the proof of \cite[Thm 1.2]{DG} that there exists $\delta_1>0$
such that if $k$ is sufficiently large, then $I - \hat{A}_{0,k}' (\overline{u}^k,
\overline{v}^k)$ is invertible and for $(u^k, v^k) \in \npair \cap C$, 
$$\mbox{index}_{\conez} (\hat{A}_{1,k}' (\uk, \vk), 0) = \mbox{index}_{\conez}(\hat{A}_{0,k}' (\overline{u}^k,
\overline{v}^k), 0)$$
since the linearisations $\hat{A}_{1,k}', \hat{A}_{0,k}'$ act in the space $\conez$, with homogeneous
Dirichlet boundary conditions, as in \cite{DG}. Since we prove that Problem (\ref{endpoint}) has a unique solution and
since the indices of non-degenerate fixed points of nonlinear maps
are equal to those of the corresponding linearised maps, it follows that
solutions $(\uk, \vk)$ of $\pks$ in $\npair$ all have the same index $\gamma$, where
$\gamma = \pm 1$. Now let $\cal N$ be the number of stationary
solutions of Problem $P_k$. We deduce from (\ref{degindexfinal}) that ${\cal N}
\times \gamma = \mbox{index}_K(B_2, w_0) = \pm 1$, which in turn implies
that ${\cal N} = 1$.\qed

\section{On the non-degeneracy condition}

Here we discuss the key non-degeneracy condition (\ref{nondegen}). Consider solutions
$w \in \wtp$ of the stationary limit problem $\stat$.
Note first that, in contrast to the case when $w$ satisfies homogeneous Neumann boundary
conditions (see a remark in \cite[p 472]{DZ}), (\ref{nondegen}) no longer always holds
in one space dimension. To see this, let $\Omega$ be an interval, say $\Omega = (0,1)$.
If $m_1 = m_2 = 0$ on $\partial \Omega$, it follows from, for example, \cite{YH}, that there
may be many solutions of $\stat$. Now suppose that $(\alpha m_1 - m_2)(0) >1$ and $(\alpha m_1 - m_2)(1) < -1$.
Then a solution $w$ of $\stat$ must be decreasing in $x \in (0,1)$. This is because the form of
$h$ (see (\ref{defh}) and {\bf (a)}) forbids a local maximum (resp. minimum)
of $w(x)$ at $x_0$ if $w(x_0)>1$ (resp. $< -1$) and the fact that $\stat$ has only even order
derivatives together with  uniqueness for initial-value problems gives that $w$ must be
symmetric about any critical point. Now define $\gamma_w (x) = w'(x)^2 + H(w(x))$, where $H$
is a primitive of $h$, and note that $\gamma_w$ is independent of $x$ for a given $w$. Hence
if $w, \hat{w}$ are two solutions of $\stat$ with $w(0) = \hat{w}(0) = (\alpha m_1 - m_2)(0) > 1$ and
$| w'(0)| \geq |\hat{w}'(0)|$, then $|w'(x)| \geq |\hat{w}'(x)|$ at any $x$ at which $w(x) = \hat{w} (x)$.
So if $(\alpha m_1 - m_2)(1) < -1$, this, together with the fact that $w, \hat{w}$
are both decreasing, yields that $w \equiv \hat{w}$, and thus $\stat$ has at most one solution in this case.
This shows that (\ref{nondegen}) cannot hold for every $m_1, m_2 \geq 0$, because if it did,
the number of solutions of $\stat$ would be preserved as $m_1, m_2$ varied, by the inverse function
theorem applied to a suitable projection (in fact, to the mapping $\pi_V|_{F^{-1}(0)}$ defined
in the proof of Theorem \ref{thmnondegen}).

However, it is possible to prove some results on non-degeneracy holding for all solutions
of $\stat$ for a generic set of boundary data. We use ideas from \cite{ST} and \cite{Dgen}.
\cite{ST} prove that if $h \in C^1$ and $h(0) =0$ then for generic $\phi$, the equation
$\Delta w + h(w) =0$ in $\Omega$, $w=\phi$ on $\partial \Omega$, has only non-degenerate solutions.
\cite{Dgen} extends their main ideas to $h$ with possible discontinuities in $h'$ in the context
of generic domain (rather than boundary data) dependence. Note that $h$ defined in (\ref{defh})
is locally Lipschitz but is not in general $C^1$.

\begin{theorem}
\label{thmnondegen}
\begin{description}
\item[(a)] There is a dense subset ${\cal A}$ of $\{ y \in \wtp: \Delta y=0 ~\mbox{in}~ \Omega \}$
such that if $\phi \in {\cal A}$ and non-negative $m_1, m_2 \in \wtp$  are such that $\alpha m_1 - m_2 =
\phi$ on $\partial \Omega$, then every solution
of $\stat$ is non-degenerate.
\item[(b)] Suppose that $\partial \Omega = \Gamma_1 \cup \Gamma_2$, where $\Gamma_i$ is closed
in $\partial \Omega$, $i=1,2$. Fix $\psi \in \wtp$ such that $\psi |_{\Gamma_1} >0$, $\psi |_{\Gamma_2}=0$
and $\Delta \psi =0$ in $\Omega$, and suppose that $m_1 = \psi$ on $\partial \Omega$. Then there is a dense
subset ${\cal B}$ of $\{ y \in \wtp: y|_{\Gamma_1} =0, ~y |_{\Gamma_2}>0, ~\Delta y =0~\mbox{in}~ \Omega \}$
such that if $\phi \in {\cal B}$ and $m_2 = \phi$ on $\partial \Omega$, then every solution of $\stat$ is
non-degenerate.
\end{description}
\end{theorem}

\pf The overall structure of the proofs is similar to that of \cite[Thm 1]{Dgen} and \cite[Thm 3.1]{ST}
and we give sketches here. The idea is to apply the version of Sard's theorem from \cite{Q} to a suitable
map. First consider {\bf (a)} and define
$$ \begin{array}{l}
X=U=W^{2,p}_0 (\Omega), \\
Y= \{ y \in \wtp : \Delta y =0 ~\mbox{in}~ \Omega \}, \\
V = Y \setminus \{ 0 \}, \\
Z = L^p (\Omega),
\end{array}
$$
(we import the notation for spaces from \cite{Dgen, ST} for ease of reference).
Define $F: U \times V \rightarrow Z$ by
$$ F(u,v) = \Delta u + h(u+v), ~~~(u,v) \in U \times V. $$
Note first that $u+v$ is not identically zero for any $(u,v) \in U \times V$. This,
together with \cite{CF}, gives that for any solution $(u,v)$ of $F(u,v)=0$, $u+v$ only
takes the value zero on a set of measure zero. It follows, arguing as in \cite[p 468-469]{DD}
and \cite[p 248]{Djde}, that $F$ is strictly differentiable at $(u,v)$ whenever $F(u,v)=0$
(see \cite[p 48]{Cartan} for the definition of strictly differentiable). This is the key property
that allows $h$ to be only locally Lipschitz - see \cite{Dgen}.

Next note two technical properties of $F$. Firstly, $F(\cdot, v):u \mapsto F(u,v)$ is Fredholm
of index zero. The mapping $u \mapsto \Delta u$ is an isomorphism from $X$ onto $Z$,
and $u \mapsto h'(u^0+v^0)u$ is a linear compact operator from $X$ into $Z$ for each $(u^0, v^0)
\in U \times V$ (note that $u^0+v^0 \neq 0 ~a.e.$, so $h'(u^0, v^0)$ makes sense, and that
$h'(u^0 + v^0) \in L^{\infty} (\Omega)$, by {\bf (a)} and (\ref{defh}) ).
Secondly, $F$ is proper, in the sense that the set of $u \in U$ such that $F(u,v)=0$ with
$v$ belonging to a compact set in $Y$ is relatively compact in $Y$. This is because if $v_n
\rightarrow v$ in $\wtp$ and $F(u_n, v_n)=0$, then since $\{ u_n+v_n \}_{n=1}^{\infty}$
is bounded in $\wtp$, $\{ u_n \}_{n=1}^{\infty}$ is bounded in $\wtp$ and has a subsequence $u_{n_k}
\rightarrow u$ in $C(\overline{\Omega})$ (since $p > N/2$). Since $h$ is continuous,
it follows that $h(u_{n_k}+v_{n_k})$ converges in $L^p(\Omega)$, so $\{\Delta( u_n + v_n) \}_{n=1}^{\infty}$
is relatively compact in $L^p(\Omega)$, and thus $\{ u_n + v_n \}_{n=1}^{\infty}$, and so
$\{ u_n \}_{n=1}^{\infty}$, is relatively compact in $W^{2,p} (\Omega)$.

We also need to check that zero is a regular value of $F$. As in \cite{Dgen, ST}, it suffices to show that if
$F(u^0, v^0)=0$ and $u \in U=W^{2,p}_0 (\Omega)$ satisfies $\Delta u + h'(u^0+v^0)u=0~~a.e.$ with
$$\int_{\Omega} h'(u^0+v^0)yu~dx =0 ~~~\mbox{for all}~~y \in Y, $$
then $u \equiv 0$. But
$$ \int_{\Omega} h'(u^0+v^0)yu~dx = - \int_{\Omega} y \Delta u~dx = - \int_{\partial \Omega} y
\frac{\partial u}{\partial \nu} ~dS~~~\mbox{for all}~~y \in Y $$
implies that $\partial u / \partial \nu =0$ on $\partial \Omega$. Since $u=0$ on $\partial \Omega$, as
$u \in W^{2,p}_0 (\Omega)$, it follows as discussed in \cite[p 144]{Dgen} that $u \equiv 0$ in
$\Omega$, as required.

Now the strict differentiability of $F$ and zero being a regular value give that $F^{-1}(0)$ is a
$C^1$-manifold (see \cite[Lem 1]{Dgen}). Also, $\pi_V|_{F^{-1}(0)} : F^{-1}(0) \rightarrow Y$ is a
$C^1$-Fredholm map of index zero,  $v^0$ is a regular
value of $\pi_V|_{F^{-1}(0)}$ if and only if $0$ is a regular value of $F(\cdot, v^0)$, and $\pi_V|_{F^{-1}(0)}$ is proper (see \cite[part {\em (i)} of the proof of
Thm 1.2 and the Appendix]{ST}). Here $\pi_V$ denotes the usual projection of $U \times V$ onto the
second factor and note that $0$ being a regular value of $F(\cdot, v^0)$ says precisely that 
every solution of $\stat$ equal to $v^0$ on $\partial \Omega$ is non-degenerate. 
The result follows by applying the version of Sard's theorem in \cite{Q} to $\pi_V|_{F^{-1}(0)}$.  (See \cite[p 144]{Dgen} for more detail on
these concluding arguments.)\\

\noindent For part {\bf (b)}, set
$$ \begin{array}{l}
X=U=W^{2,p}_0(\Omega), \\
Y= \{ w \in \wtp: w|_{\Gamma_1} =0, ~~\Delta w =0~~\mbox{in}~~\Omega \},\\
V= \{ w \in Y: w|_{\Gamma_2} >0 \}, \\
Z= L^p(\Omega),
\end{array}
$$
and define $F : U \times V \rightarrow Z$ by
$$ F(u,v) = \Delta u + h(u+v+\psi), ~~~(u,v) \in U \times V. $$
Most of the properties noted for {\bf (a)} follow in the same way here. To see that
zero is a regular value of $F$, note that if $F(u^0, v^0)=0$, $u \in U=W^{2,p}_0(\Omega)$
satisfies $ \Delta u + h'(u^0+v^0+\psi)u =0~~~a.e$ and
$$ \int_{\Omega}h'(u^0+v^0+\psi)yu~dx =0 ~~~\mbox{for all}~~y \in Y,$$
then
\begin{equation}
\label{last}
 \int_{\Omega} h'(u^0+v^0+\psi)yu~dx = - \int_{\Gamma_2} y \frac{\partial u}{\partial \nu}~dS =0
~~~\mbox{for all}~~y \in Y,
\end{equation}
which yields that $\partial u / \partial \nu =0$ on $\Gamma_2$. This is enough to deduce that $u \equiv 0$
in $\Omega$ - again, see \cite[p 144]{Dgen}. The remainder of the proof is the same as for {\bf (a)}. \qed

\noindent {\bf Remarks.}  {\em (i)} Note that if $N=1$, the boundary data is necessarily constant on each component of $\partial \Omega$. In particular, {\bf (b)} implies that if $\Omega = (0,1)$ and we fix
$m_1 (0) >0$, $m_1(1)=0$ and $m_2(0)=0$, then there is a dense set of positive constant values
for $m_2(1)$ for which every solution of $\stat$ is non-degenerate. When $N>1$ and $0,1$ are
replaced by $\Gamma_1, \Gamma_2$ respectively, our method does not give a corresponding result for boundary data constant on each of the components $\Gamma_1, \Gamma_2$ of $\partial \Omega$
since in that case we can no longer deduce from (\ref{last}) that $\partial u / \partial \nu  \equiv 0$
on $\Gamma_2$.\\[5pt]
{\em (ii)} Application of the inverse function theorem to $\pi_V|_{F^{-1}(0)}$ at a non-degenerate
solution $u^0 + v^0$ of $\stat$ yields isolatedness in $W^{2,p}(\Omega)$ of solutions of
$\stat$ equal to $v^0$ on $\partial \Omega$, from which isolatedness in $L^2(\Omega)$
follows (because  solutions of $\stat$ belong to a bounded set in $L^{\infty}(\Omega)$ 
for a given $v^0$, so non-isolatedness in $L^2(\Omega)$ would imply
non-isolatedness in $L^p(\Omega)$, and since solutions of $\stat$ satisfy (\ref{stationary})
and $h$ is locally Lipschitz, 
non-isolatedness in $L^p(\Omega)$ would imply non-isolatedness in $W^{2,p}(\Omega)$).

\bigskip \bigskip \bigskip
\noindent
{\sc E.C.M. Crooks} {\em (corresponding author)}\\
Mathematical Institute, 24-29 St Giles', Oxford, OX1 3LB, U.K.\\
\texttt{crooks@maths.ox.ac.uk}\\ \\
{\sc E.N. Dancer}\\
School of Mathematics and Statistics,  University of Sydney, NSW 2006, Australia.\\
\texttt{normd@maths.usyd.edu.au}\\ \\
{\sc D. Hilhorst}\\
CNRS and Laboratoire de Math{\'{e}}matiques,  \\Universit{\'{e}}
de Paris-Sud, (b\^{a}t. 425), 91405 Orsay-Cedex, France. \\
\texttt{Danielle.Hilhorst@math.u-psud.fr}

\end{document}